\documentclass[12pt]{article}
\usepackage{amsfonts,latexsym,graphicx,
amssymb,makeidx,amsthm,amsmath}
\newcommand{\bO}{\mbox{\boldmath $O$}}
\newcommand{\bH}{\mbox{\boldmath $H$}}
\newcommand{\bC}{\mbox{\boldmath $C$}}
\newcommand{\bR}{\mbox{\boldmath $R$}}

\newcounter{f1}
\newcounter{f2}
\newcounter{f3}
\newcounter{f4}
\newcounter{f5}
\newcounter{f6}
\newcounter{f7}
\newcounter{f8}
\newcounter{f9}
\newcounter{f10}
\newcounter{f11}
\newcounter{f12}
\renewcommand{\refname}{Reference}
\title{
Orbit decomposition of Jordan 
matrix algebras of order three 
under the automorphism groups 
}
\author{Akihiro Nishio and Osami Yasukura}
\date{}
\begin{document}

\maketitle

{\bf Abstract}. 
The orbit decomposition is given 
under the automorphism group on 
the real split Jordan algebra of 
all hermitian matrices of order three 
corresponding to any real split composition algebra, 
or the automorphism group on the complexification, 
explicitly, in terms of 
the cross product of H.~Freudenthal 
and the characteristic polynomial.

\bigskip

\begin{center}
{\bf 0. Introduction.} 
\end{center}

Let $\mathcal{J}'$ 
be a split exceptional simple Jordan algebra 
over a field $\mathbb{F}$ of characteristic not two, 
that is, 
the set of all hermitian matrices of order three 
whose elements are split octonions over $\mathbb{F}$ 
with the Jordan product.  
And let $G'$ be the automorphism group of $\mathcal{J}'$. 
N.~Jacobson \cite[p.389, Theorem 10]{Jn1968} 
found that 
$X, Y \in \mathcal{J}'$ 
are in the same $G'$-orbit if and only if 
$X, Y$ admit the same minimal polynomial 
and the same generic minimal polynomial, 
by imbedding a generating subalgebra 
with the identity element $E$ 
in terms of the Jordan product 
into a special Jordan algebra. 
When $\mathbb{F} = \mathbb{R}$, 
the field of all real numbers, 
some elements of ${\cal J}'$ 
are not diagonalizable under 
the action of $G' = F_{4(4)}$,  
since ${\cal J}'$ admits a $G'$-invariant 
non-defnite $\mathbb{R}$-bilinear form
such that the restriction to the subspace of 
all diagonal elements is positive-definite 
\cite[Theorem 2]{MY2001}, 
although every element of ${\cal J}'$ 
is diagonalizable under the action of 
a linear group $E_{6(6)}$ containing $F_{4(4)}$ 
on ${\cal J}'$ 
by \cite{Jn1961} (cf. \cite{Ks2002})  
or under the action of the maximal compact subgroup 
$Sp(4)/\mathbb{Z}_2$ of $E_{6(6)}$ on ${\cal J}'$ 
given by \cite{SY1979}.

This paper presents a concrete orbit 
decomposition under the automorphism group 
on a real split Jordan algebra 
of all hermitian matrices of order three 
corresponding to any real split composition algebra, 
or the complexification of it, 
that is special or exceptional as a Jordan algebra. 
As a result, 
$X, Y \in \mathcal{J}'$ 
are in the same $G'$-orbit if and only if 
$X, Y$ admit the same dimension of the generating 
subspace with $E$ by the cross product 
\cite{Fh1953} 
and the same characteristic polynomial, 
which gives a simplification for 
N.~Jacobson \cite{Jn1968}'s polynomial invariants 
on $G'$-orbits when $\mathbb{F} = \mathbb{R}$ 
or the field of all complex numbers 
$\mathbb{C}$. 
To state the main results more precisely, 
let us give the precise notations:

Put 
$\mathbb{F} := \mathbb{R}$ 
or $\mathbb{C}$. 
Let $V$ be an $\mathbb{F}$-linear space, 
and $\mathrm{End}_{\mathbb{F}}(V)$ 
(or $\mathrm{GL}_{\mathbb{F}}(V)$) 
denote the set of all 
$\mathbb{F}$-linear endomorphisms 
(resp. automorphims) on $V$. 
For a mapping 
$f: V \rightarrow V$ 
and $c \in \mathbb{F}$, 
put $V_{f, c} := \{ v \in V |~f(v)= c v\}$ 
and $V_{f, 1} := V_f$.  
For a subgroup $G$ of $\mathrm{GL}_{\mathbb{F}}(V)$, 
let $G^{\circ}$ 
be the identity connected component of $G$. 
For $v\in V$ and a mapping 
$\phi: V \rightarrow V$, 
put 
$\mathcal{O}_G(v) 
:= \{ \alpha (v) |~\alpha \in G \}$, 
$G_v := \{ \alpha \in G |~\alpha (v) = v \}$ 
and $G^{\phi} := 
\{ \alpha \in G 
|~\phi \circ \alpha = \alpha \circ \phi \}$. 
For a subset $W$ of $V$, put 
$G_W := \{ \alpha \in G |~\{ \alpha w |~w \in W \} = W \}$. 
For positive integers $n, m$, let $M(n, m; V)$ 
be the set of all $n \times m$-matrices 
with entries in $V$. 
Put $V^n := M(n, 1; V), V_m := M(1, m; V)$ 
and $M_n( V) := M(n, n; V)$. 
Since $V$ can be considered as an $\mathbb{R}$-linear space, 
the complexification is defined as 
$V^{\mathbb{C}} := V \otimes_{\mathbb{R}} \mathbb{C} 
=  V \oplus \sqrt{-1} V$ 
with an $\mathbb{R}$-linear conjugation: 
$\tau: V^{\mathbb{C}} \rightarrow V^{\mathbb{C}}; 
v_1 + \sqrt{-1} v_2 
\mapsto v_1 - \sqrt{-1} v_2$ 
$(v_1, v_2 \in V)$. 
For any 
$\alpha \in \mathrm{End}_{\mathbb{R}}(V)$, 
put 
$\alpha^{\mathbb{C}}: 
V^{\mathbb{C}} \longrightarrow V^{\mathbb{C}}; 
v_1 + \sqrt{-1} v_2 
\mapsto (\alpha v_1) + \sqrt{-1} (\alpha v_2)$ 
such that 
$\alpha^{\mathbb{C}} \tau = \tau \alpha^{\mathbb{C}}$, 
which is identified with 
$\alpha \in \mathrm{End}_{\mathbb{R}}(V)$: 
$\alpha = \alpha^{\mathbb{C}}$.

By W.R.~Hamilton, 
{\it the quaternions} 
is defined as an $\mathbb{R}$-algebra 
$\mathbb{H} := \oplus_{i=0}^3 \mathbb{R} e_i$ 
given as 
$e_0 e_i = e_i e_0 = e_i,~e_i^2 = - e_0$ 
($i \in \{ 1, 2, 3 \}$); 
$e_k e_{k+1} = - e_{k+1} e_k = e_{k+2}$ 
(where $k, k+1, k+2 \in \{ 1, 2, 3 \}$ 
are counted modulo $3$) 
with the unit element 
$1 := e_0$ 
and the conjugation 
$\overline{\sum_{i=0}^3 x_i e_i} 
= x_0 e_0 - \sum_{k = 1}^3 x_k e_k$, 
which contains {\it the complex numbers} 
$\mathbb{C} := \mathbb{R} e_0 \oplus \mathbb{R} e_1$ 
and the real numbers 
$\mathbb{R} := \mathbb{R} e_0$ 
as $\mathbb{R}$-subalgebras. 
By A.~Cayley and J.T.~Graves, 
{\it the octanions} is defined as 
a non-associative $\mathbb{R}$-algebra 
$\bO := \mathbb{H} \oplus \mathbb{H} e_4$ 
given as follows \cite{Dle1919}:  

\[
(x \oplus y e_4) (x' \oplus y' e_4) 
:= (x x' - \overline{y'} y) 
\oplus 
(y \overline{x'} + y' x) e_4 
\]

\noindent 
with the $\mathbb{R}$-linear basis 
$\{ e_i |~i = 0, 1, 2, 3, 4, 5, 6, 7 \}$, 
where the numbering is given as 
$e_5 := e_1 e_4, e_6 := - e_2 e_4, 
e_7 := e_3 e_4$ 
after \cite[p.127]{Yi1967}, 
\cite[p.20]{Dd1978} or \cite{Ms1989}. 
Put $\bH := \mathbb{C} \oplus \mathbb{C} e_4$ 
and $\bC := \mathbb{R} \oplus \mathbb{R} e_4$. 
For $K := \bO, \bH, \bC, \mathbb{R}$, 
put $d_{K} := \mathrm{dim}_{\mathbb{R}} K$. 
And put 
$\sqrt{-1} := e_0 \otimes e_1 
\in K^{\mathbb{C}} := K \otimes_{\mathbb{R}} \mathbb{C}$ 
with the identification 
$K = K \otimes e_0 \subset K^{\mathbb{C}}$. 
Then 
$K^{\mathbb{C}} =  K \oplus \sqrt{-1} K$ 
is {\it split} (i.e. non-division) 
as a ${\mathbb{C}}$-algebra with 
$\tau: K^{\mathbb{C}} \rightarrow K^{\mathbb{C}}; 
x + \sqrt{-1} y \mapsto x - \sqrt{-1} y
~(x, y \in K)$ 
as the complex conjugation with respect to the real form $K$. 
Put

\begin{eqnarray*}
& &\gamma: 
\bO^{\mathbb{C}} \longrightarrow \bO^{\mathbb{C}}; 
\sum_{i=0}^{7} x_i e_i \mapsto 
\sum_{i=0}^{3} x_i e_i - \sum_{i=4}^{7} x_i e_i;~{\rm and} \\
& &\epsilon: 
\bO^{\mathbb{C}} \longrightarrow \bO^{\mathbb{C}}; 
x := \sum_{i=0}^{7} x_i e_i 
\mapsto \bar{x} := x_0 - \sum_{i=1}^{7} x_i e_i
\end{eqnarray*}

\noindent 
as $\mathbb{C}$-linear conjugations with respect to 
$\mathbb{H}^{\mathbb{C}}$ 
and $\mathbb{R}^{\mathbb{C}}$, 
respectively. 
And a ${\mathbb{C}}$-bilinear form are defined on $\bO^{\mathbb{C}}$ as 
$(x | y) := 
(x \bar{y} + \overline{x \bar{y}}) / 2 
= \sum_{i=0}^{7} x_i y_i \in {\mathbb{C}}$. 
The restrictions of $\gamma, \epsilon$ 
and $(x | y)$ on $K^{\mathbb{C}}$ 
are also well-defined and denoted by the same letters. 
Then $K^{\mathbb{C}}$ is a composition ${\mathbb{C}}$-algebra 
with respect to the norm form given by  
$N(x) := (x | x)$ \cite[\S \Roman{f1}.3]{Dd1978}, 
because of 
$N((x \oplus y e_4) (x' \oplus y' e_4)) 
- N(x \oplus y e_4) N(x' \oplus y' e_4) 
= 2 \{ (y \overline{x'} | y' x) - (x x' | \overline{y'} y) \} 
= 2 (\overline{y'} (y \overline{x'}) - (\overline{y'} y) \overline{x'} |~x) 
= 0$ 
since $\mathbb{H}^{\mathbb{C}}$ 
is an associative composition algebra with respect to 
$N$ \cite[\S 6.4]{CS2003}. 
And 
$K = (K^{\mathbb{C}})_{\tau}$ 
is a division composition 
$\mathbb{R}$-algebra with the norm form 
$N(x)$ such that $a^{-1} = \bar{a}/N(a)$ 
for $a \ne 0$. 

Put 
$K' := (K^{\mathbb{C}})_{\tau \gamma}$ 
as a composition $\mathbb{R}$-algebra 
with the norm form $N(x)$ such that 
$(K')_{\gamma} = K_{\gamma} 
= (K')_{\tau} = K' \cap K$. 
Precisely, 
$\bO' = \{ \sum_{i=0}^3 x_i e_i 
+ \sum_{i=4}^7 x_i \sqrt{-1} e_i 
|~x_i\in \mathbb{R} \}$ 
is the $\mathbb{R}$-algebra 
of {\it the split-octanions} 
containing the $\mathbb{R}$-subalgebra 
$\bH' 
= \{ \sum_{i=0}^1 x_i e_i 
+ \sum_{i=4}^5 x_i \sqrt{-1}e_i 
|~x_i \in \mathbb{R} \}$ 
of {\it the split-quaternions} and 
the $\mathbb{R}$-subalgebra 
$\bC' = \{ x_0 + x_4 \sqrt{-1} e_4 
|~x_i \in \mathbb{R} \}$ 
of {\it the split-complex numbers} 
such that 
$\bO' \cap \bO = \mathbb{H}$, 
$\bH' \cap \bH = \mathbb{C}$ 
and $\bC' \cap \bC = \mathbb{R}$. 
Then 
${K'}^{\mathbb{C}} = K' \oplus \sqrt{-1} K' 
= K^{\mathbb{C}}$ 
as a $\mathbb{C}$-subalgebra of $\bO^{\mathbb{C}}$.

Put 
$\tilde{K} := K, K'$ 
(or 
${K'}^{\mathbb{C}}$, 
$K^{\mathbb{C}}$)  
with $\mathbb{F} := \mathbb{R}$ (resp. ${\mathbb{C}}$) 
and $d_{\tilde{K}} := \mathrm{dim}_{\mathbb{F}} \tilde{K}$. 
For $A \in M_n( \tilde{K})$ with the $(i, j)$-entry 
$a_{i j} \in \tilde{K}$, let 
$^t A, \tau A, \epsilon A 
\in M_n( \tilde{K})$ 
be the transposed, $\tau$-conjugate, 
$\epsilon$-conjugate matrix of $A$ 
such that the $(i, j)$-entry is equal to 
$a_{j i}, \tau (a_{i j}), \epsilon (a_{i j})$, 
respectively, with the trace 
$\mathrm{tr}(A) 
:= \sum_{i=1}^n a_{i i} \in \mathbb{F}$, 
and the adjoint matrix 
$A^*:=~^t (\epsilon A) \in M_n( \tilde{K})$. 
Let denote the set of 
all hermitian matrices of order three 
corresponding to $\tilde{K}$ as follows: 

\[
\mathcal{J}_3(\tilde{K}) 
:= \{ X \in M_3(\tilde{K}) |~X^* = X \}  
\]

\noindent 
with an $\mathbb{F}$-bilinear Jordan algebraic product 
$X \circ Y := \frac{1}{2}(XY + YX)$, 
the identity element 
$E := \mathrm{diag}(1,1,1)$ 
and an $\mathbb{F}$-bilinear symmetric form 
$(X | Y) := \mathrm{tr}(X \circ Y) \in \mathbb{F}$. 
After H.~Freudenthal \cite{Fh1953} 
(cf. \cite[(7.5.1)]{Fh1951}, \cite{Yi1959}, 
\cite{Jn1960}, \cite[p.232, (47)]{Jn1968}, 
\cite{Yi1975}), 
{\it the cross product} on 
$\mathcal{J}_3(\tilde{K})$ is defined 
as follows:

\[
X \times Y := 
X \circ Y 
- \frac{1}{2} ( 
\mathrm{tr}(X) Y + \mathrm{tr}(Y) X 
- ( \mathrm{tr}(X) \mathrm{tr}(Y) - (X | Y) ) E ) 
\]

\noindent 
with $X^{\times 2} := X \times X$ 
as well as an $\mathbb{F}$-trilinear form 
$(X | Y | Z) := (X \times Y | Z)$ 
and the determinant 
$\mathrm{det}(X) 
:= \frac{1}{3}(X | X | X) \in \mathbb{F}$ 
on $\mathcal{J}_3(\tilde{K})$ 
(cf. \cite[p.163]{Fh1964}). 
Put 
$E_i := \mathrm{diag}(\delta_{i1}, \delta_{i2}, \delta_{i3})$ 
for $i \in \{ 1, 2, 3 \}$ 
with the Kronecker's delta $\delta_{ij}$. 
For $x \in \tilde{K}$, put

\[
F_1(x) := 
\begin{pmatrix}
   0& 0 & 0\\
   0 & 0 & x\\
   0& \overline{x} & 0
\end{pmatrix},
~F_2(x) := 
\begin{pmatrix}
   0& 0 & \overline{x}\\
   0 & 0 & 0\\
   x& 0 & 0
\end{pmatrix},
~F_3(x) := 
\begin{pmatrix}
   0& x & 0\\
   \overline{x} & 0 & 0\\
   0& 0 & 0\\
\end{pmatrix}. 
\]

\noindent 
For 
$x \in K^{\mathbb{C}} 
= \mathbb{R}^{\mathbb{C}}, 
\bC^{\mathbb{C}}, 
\bH^{\mathbb{C}}$ 
or 
$\bO^{\mathbb{C}}$, 
put 
$M_1(x), M_{2 3}(x) 
\in \mathcal{J}_3(K^{\mathbb{C}})$ 
such as

\[
M_1(x) 
:= (x | 1) (E_2 - E_3) 
+ F_1(\sqrt{-1} x),
~~M_{2 3}(x) 
:= F_2(\sqrt{-1} \overline{x}) + F_3(x) 
\]

\noindent 
with 
$M_1 := M_1(1)$, 
$M_{2 3} := M_{2 3}(1)$. 
For $x \in K' = \bC', \bH'$ or $\bO'$, 
put 
$M_{1'}(x), M_{2' 3}(x) 
\in \mathcal{J}_3(K')$ 
such as 

\[
M_{1'}(x) := (x | 1) (E_2 - E_3) 
+ F_1(\sqrt{-1} e_4 x),
~~M_{2' 3}(x) := 
F_2(- \sqrt{-1} e_4 \overline{x}) + F_3(x) 
\]

\noindent 
with 
$M_{1'} := M_{1'}(1)$, 
$M_{2' 3} := M_{2' 3}(1)$. 
For $x \in \tilde{K}$, 
let denote 

\begin{eqnarray*}
& &\tilde{M}_{1}(x) := M_1(x) 
~({\rm when}~\tilde{K} = K^{\mathbb{C}})
~{\rm or}~M_{1'}(x) 
~({\rm when}~\tilde{K} = K'), \\
& &\tilde{M}_{2 3}(x) := M_{2 3}(x) 
~({\rm when}~\tilde{K} = K^{\mathbb{C}})
~{\rm or}~M_{2' 3}(x) 
~({\rm when}~\tilde{K} = K'); \\
& &\tilde{M}_1 := M_1 
~({\rm when}~\tilde{K} = K^{\mathbb{C}})
~{\rm or}~M_{1'}
~({\rm when}~\tilde{K} = K'), \\
& &\tilde{M}_{2 3} := M_{2 3} 
~({\rm when}~\tilde{K} = K^{\mathbb{C}}) 
~{\rm or}~M_{2' 3}
~({\rm when}~\tilde{K} = K'). 
\end{eqnarray*}

\noindent
And denote 

\begin{eqnarray*}
& &
\mathcal{P}_2(\tilde{K}) 
:= \{ X \in \mathcal{J}_3(\tilde{K}) 
|~X^{\times 2} = 0,~\mathrm{tr}(X) = 1 \}, \\
& &\mathcal{J}_3(\tilde{K})_0 
:= \{ X \in \mathcal{J}_3(\tilde{K}) 
|~\mathrm{tr}(X) = 0 \}, \\
& &\mathcal{M}_1(\tilde{K}) 
:= \{ X \in \mathcal{J}_3(\tilde{K})_0 
|~X \ne 0,~X^{\times 2} = 0 \}, \\
& &\mathcal{M}_{23}(\tilde{K}) 
:= \{ X \in \mathcal{J}_3(\tilde{K})_0 
|~X^{\times 2} \ne 0,
~\mathrm{tr}(X^{\times 2}) 
= \mathrm{det}(X) = 0 \}. 
\end{eqnarray*}

\noindent 
When $\tilde{K} = K$, 
$\mathcal{P}_2(\tilde{K})$ 
has a structure of Moufang projective plane  
\cite[p.162, 4.6, 4.7]{Fh1964}, 
the algebraization method of which motivates 
to define the cross product on 
$\mathcal{J}_3(\tilde{K})$ 
for any $\tilde{K}$. 
The automorphism group of 
$\mathcal{J}_3(\tilde{K})$ 
with respect to the 
$\mathbb{F}$-bilinear 
Jordan product $X \circ Y$ 
is denoted as follows: 

\[
G(\tilde{K}) 
:= \mathrm{Aut}(\mathcal{J}_3(\tilde{K})) 
= \{ \alpha \in 
GL_{\mathbb{F}}(\mathcal{J}_3(\tilde{K})) 
|~\alpha(X \circ Y) 
= \alpha X \circ \alpha Y \},  
\]

\noindent 
which is a complex 
(resp. compact; real split) 
simple Lie group of type $(F_4)$ 
(resp. $(F_{4(-52)})$; $(F_{4(4)})$) 
when 
$\tilde{K} 
= {\bO'}^{{\mathbb{C}}} = \bO^{{\mathbb{C}}}$ 
(resp. $\bO$; $\bO'$) 
by C.~Chevalley 
and R.D.~Schafer \cite{CS1950} 
(resp. \cite{Fh1951}, 
\cite[p.161]{Fh1964}, 
\cite[p.206, (2), (3)]{Ms1989}; 
\cite{Yi1977}). 
When $K = \mathbb{R}, \bC$ or $\bH$, 
the group 
$G(\tilde{K})$ 
is a simple Lie group of type 
$(A_1), (A_2)$ or $(C_3)$, 
respectively 
(cf. \cite[p.165]{Fh1964}). 
Put 
$\gamma: 
\mathcal{J}_3(\tilde{K}) 
\longrightarrow \mathcal{J}_3(\tilde{K}); 
X \mapsto \gamma X$ 
such that 
$\gamma X := 
\sum_{i = 1}^3 (\xi_i E_i + F_i(\gamma x_i))$ 
for 
$X = \sum_{i = 1}^3 (\xi_i E_i + F_i(x_i)) 
\in \mathcal{J}_3(\tilde{K})$. 
Put 
$\tau: \mathcal{J}_3(\tilde{K}) 
\longrightarrow \mathcal{J}_3(\tilde{K}); 
X \mapsto \tau X$ 
such as 
$\tau X := \sum_{i=1}^3 ((\tau \xi_i) E_i + F_i(\tau x_i))$. 
Then 
$\tau \in GL_{\mathbb{R}}(\mathcal{J}_3(\tilde{K}))$ 
such that 
$\tau (X \circ Y) = (\tau X) \circ (\tau Y), 
\tau (X \times Y) = (\tau X) \times (\tau Y)$,
$\mathrm{tr}(\tau X) = \tau (\mathrm{tr}(X))$, 
$(\tau X | \tau Y) = \tau (X | Y)$ 
and 
$\mathrm{det}(\tau X) = \tau (\mathrm{det} X)$, 
and that 
$\tau^2 = \mathrm{id}$, 
$\mathcal{J}_3(\tilde{K}) 
= \mathcal{J}_3(\tilde{K})_{\tau} \oplus 
\mathcal{J}_3(\tilde{K})_{-\tau}$ 
and 
$\mathcal{J}_3(K^{\mathbb{C}})_{-\tau} 
= \sqrt{-1} \mathcal{J}_3(K)$, 
so that 
$G(K) \equiv \{ \alpha^{\mathbb{C}} |~\alpha \in G(K) \} 
= G(K^{\mathbb{C}})^{\tau}$.

For $X \in \mathcal{J}_3(\tilde{K})$ 
and the indeterminate $\lambda$, 
put $\varphi_X (\lambda) := \lambda E - X$. 
Then {\it the characteristic polynomial of $X$} 
is defined as the polynomial  
$\Phi _X(\lambda) := 
\mathrm{det}(\varphi_X (\lambda))$ 
of $\lambda$ with degree 3 and the derivative 
$\Phi_X'(\lambda)$ 
is 
$\frac{d}{d \lambda} \Phi_X(\lambda)$, 
so that 
$\Phi _{X}(\lambda) \equiv 
(\lambda - \lambda_1)(\lambda - \lambda_2) 
(\lambda - \lambda_3)$ 
with some 
$\lambda_1, \lambda_2, \lambda_3 \in {\mathbb{C}}$. 
In this case, the set 
$\{ \lambda_1, \lambda_2, \lambda_3 \}$ 
is said to be 
{\it the characteristic roots of $X$}. 
Put 
$\Lambda_X := 
\{ \lambda_1, \lambda_2, \lambda_3 \} 
\subset {\mathbb{C}}$ 
with 
$\# \Lambda_X \in \{ 1, 2, 3 \}$ 
and 
$V_X := 
\{ a X^{\times 2} + b X + c E 
|~a, b, c \in \mathbb{F} \}$ 
with 
$v_X := \mathrm{dim} V_X \in \{ 1, 2. 3 \}$. 

\bigskip

{\sc Proposition 0.1.} 
{\it 
Let $\tilde{K}$ be $K, K'$ or $K^{\mathbb{C}}$ 
with $K = \mathbb{R}, \bC, \bH$ or $\bO$. 
}

\medskip

\noindent
(1) 
{\it 
$G(\tilde{K}) \subseteqq 
\{ \alpha 
\in \mathrm{GL}_{\mathbb{F}}(\mathcal{J}_3(\tilde{K})) 
|~\mathrm{tr}(\alpha X) = \mathrm{tr}(X),~\alpha E = E \}$. 
And  
} 

\begin{eqnarray*} 
G(\tilde{K}) 
&=& \{ \alpha 
\in \mathrm{GL}_{\mathbb{F}}(\mathcal{J}_3(\tilde{K})) 
|~\mathrm{det} (\alpha X) 
= \mathrm{det}(X),~\alpha E = E \} \\
&=& \{ \alpha 
\in \mathrm{GL}_{\mathbb{F}}(\mathcal{J}_3(\tilde{K})) 
|~\Phi_{\alpha X}(\lambda) = \Phi_X(\lambda) \} \\
&=& \{ \alpha 
\in \mathrm{GL}_{\mathbb{F}}(\mathcal{J}_3(\tilde{K})) 
|~\mathrm{det} (\alpha X) = \mathrm{det}(X),
~(\alpha X | \alpha Y) = (X | Y) \} \\ 
&=& \{ \alpha 
\in \mathrm{GL}_{\mathbb{F}}(\mathcal{J}_3(\tilde{K})) 
|~\alpha (X \times Y) = (\alpha X) \times (\alpha Y) \}. 
\end{eqnarray*}

\noindent 
{\it 
Especially, 
$\Lambda_{\alpha X} = \Lambda_X$ 
and $v_{\alpha X} = v_X$ 
for all $X \in \mathcal{J}_3(\tilde{K})$ 
and $\alpha \in G(\tilde{K})$. 
}

\medskip

\noindent 
(2) 
{\it 
$G(\tilde{K})^{\tau}$ is a maximal compact 
subgroup of $G(\tilde{K})$. 
And 
$\gamma \in G(\tilde{K})^{\tau}_{E_1, E_2, E_3}$.
}

\medskip

\noindent 
(3) 
$\mathcal{P}_2(\tilde{K}) 
= \mathcal{O}_{G(\tilde{K})^{\circ}}(E_1)$, 

\medskip

\noindent 
(4) 
{\it 
Assume that $\tilde{K} \ne K$, i.e., 
$\tilde{K} = \bR^{\mathbb{C}}, 
\bC^{\mathbb{C}}, \bH^{\mathbb{C}}, 
\bO^{\mathbb{C}}; \bC', \bH'$ or $\bO'$. 
Then: 
}

\smallskip

(\roman{f1}) 
$\mathcal{M}_1(\tilde{K}) 
= \mathcal{O}_{G(\tilde{K})^{\circ}}(\tilde{M}_1)$, 

(\roman{f2}) 
$\mathcal{M}_{2 3}(\tilde{K}) 
= \mathcal{O}_{G(\tilde{K})^{\circ}}(\tilde{M}_{2 3})$.

\bigskip

{\sc Theorem 0.2.} 
{\it 
Let $K^{\mathbb{C}}$ 
be $\mathbb{R}^{\mathbb{C}}, 
\bC^{\mathbb{C}}, \bH^{\mathbb{C}}$ 
or $\bO^{\mathbb{C}}$. 
Then the orbit decomposition of 
$\mathcal{J}_3(K^{\mathbb{C}})$ 
over $G(K^{\mathbb{C}})$ 
or $G(K^{\mathbb{C}})^{\circ}$ 
is given as follows: 
} 

\smallskip 

\noindent 
(1) 
{\it 
Take $X \in \mathcal{J}_3(K^{\mathbb{C}})$. 
Then $\# \Lambda_X = 3, 2$ or $1$. 
}

\smallskip

\noindent 
(\roman{f1}) 
{\it 
Assume that 
$\# \Lambda_X = 3$ 
with 
$\Lambda_X = \{ \lambda_1, \lambda_2, \lambda_3 \}$. 
Then 
$\mathrm{diag}(\lambda_1, \lambda_2, \lambda_3) 
\in {\cal O}_{G(K^{\mathbb{C}})^{\circ}}(X)$ 
with $v_X = 3$.  
}

\smallskip

\noindent 
(\roman{f2}) 
{\it 
Assume that $\#\Lambda_X = 2$ with 
$\Lambda_X = \{ \lambda_1, \lambda_2 \}$ 
such that 
$\Phi_X'(\lambda_2) = 0$. 
Then $v_X = 2$ or $3$. 
Moreover:} 

\smallskip

(\roman{f2}-1) 
{\it 
$v_X = 2$ iff 
$\mathrm{diag}(\lambda_1, \lambda_2, \lambda_2) 
\in {\cal O}_{G(K^{\mathbb{C}})^{\circ}}(X)$; 
and 
}

\smallskip

(\roman{f2}-2) 
{\it 
$v_X = 3$ iff 
$\mathrm{diag}(\lambda_1,\lambda_2,\lambda_2) 
+ M_1 \in {\cal O}_{G(K^{\mathbb{C}})^{\circ}}(X)$.}

\smallskip

\noindent 
(\roman{f3}) 
{\it 
Assume that $\# \Lambda_X = 1$ 
with $\Lambda_X = \{ \lambda_1 \}$. 
Then:} 

\smallskip

(\roman{f3}-1) 
{\it 
$v_X = 1$ iff 
$\lambda_1 E \in {\cal O}_{G(K^{\mathbb{C}})^{\circ}}(X)$; 
}

\smallskip

(\roman{f3}-2) 
{\it 
$v_X = 2$ iff 
$\lambda_1 E + M_1 \in 
{\cal O}_{G(K^{\mathbb{C}})^{\circ}}(X)$; and 
}

\smallskip

(\roman{f3}-3) 
{\it 
$v_X = 3$ iff 
$\lambda_1 E + M_{2 3} \in 
{\cal O}_{G(K^{\mathbb{C}})^{\circ}}(X)$. 
}

\medskip

\noindent 
(2) 
{\it 
For $X, Y \in \mathcal{J}_3(K^{\mathbb{C}})$, 
${\cal O}_{G(K^{\mathbb{C}})^{\circ}}(X) 
= {\cal O}_{G(K^{\mathbb{C}})^{\circ}}(Y)$ 
iff $\Lambda_X = \Lambda_Y$ and $v_X = v_Y$. 
For any $X \in \mathcal{J}_3(K^{\mathbb{C}})$, 
${\cal O}(X) := {\cal O}_{G(K^{\mathbb{C}})^{\circ}}(X) 
= {\cal O}_{G(K^{\mathbb{C}})}(X)$ and 
${\cal O}(X) \cap \mathcal{J}_3(\mathbb{R}^{\mathbb{C}}) 
\ne \emptyset$. 
}

\bigskip

{\sc Theorem 0.3.} 
{\it 
Let $K'$ be $\bC', \bH'$ or $\bO'$. 
Then the orbit decomposition of 
$\mathcal{J}_3(K')$ over $G(K')$ 
or $G(K')^{\circ}$ 
is given as follows:} 

\smallskip 

\noindent 
(1) 
{\it 
Take 
$X \in \mathcal{J}_3(K')$. 
Then $\# \Lambda_X = 3, 2$ or $1$. 
}

\smallskip

\noindent
(\roman{f1}) 
{\it 
Assume that $\# \Lambda_X = 3$. 
Then $v_X = 3$. 
And 
$\Lambda_X = \{ \lambda_1, \lambda_2, \lambda_3 \}$ 
for some $\lambda_1 \in \mathbb{R}$ 
and $\lambda_2, \lambda_3 \in {\mathbb{C}}$ 
such that 
$\Lambda_X \subset \mathbb{R}$ 
or $\{ \lambda_2, \lambda_3 \} 
= \{ p \pm q \sqrt{-1} \}$ 
with some $p \in \mathbb{R}$ 
and $q \in \mathbb{R} \backslash \{ 0 \}$. 
Moreover:}

\smallskip

(\roman{f1}-1) 
{\it 
If $\Lambda_X \subset \mathbb{R}$ 
with $\lambda_1 > \lambda_2 > \lambda_3$, 
then 
$\mathrm{diag}(\lambda_1, \lambda_2, \lambda_3) 
\in {\cal O}_{G(K')^{\circ}}(X)$; and
}

\smallskip

(\roman{f1}-2) 
{\it 
If $\{ \lambda_2, \lambda_3 \} = \{ p \pm q \sqrt{-1} \}$ 
with some $p, q \in \mathbb{R}$ such that $q > 0$, then 
$\mathrm{diag}(\lambda_1, p, p) 
+ F_1(q \sqrt{-1} e_4) 
\in {\cal O}_{G(K')^{\circ}}(X)$. 
}

\smallskip

\noindent 
(\roman{f2}) 
{\it 
Assume that $\# \Lambda_X = 2$ 
with $\Lambda_X = \{ \lambda_1, \lambda_2 \}$ 
such that $\Phi_X'(\lambda_2) = 0$. 
Then $\lambda_1, \lambda_2 \in \mathbb{R}$ 
and $v_X = 2$ or $3$. 
Moreover: 
}

\smallskip

(\roman{f2}-1) 
{\it 
$v_X = 2$ iff 
$\mathrm{diag}(\lambda_1,\lambda_2,\lambda_2) 
\in {\cal O}_{G(K')^{\circ}}(X)$; 
and
}

\smallskip

(\roman{f2}-2) 
{\it 
$v_X = 3$ iff 
$\mathrm{diag}(\lambda_1, \lambda_2, \lambda_2) 
+ M_{1'} \in {\cal O}_{G(K')^{\circ}}(X)$. 
} 

\smallskip 

\noindent
(\roman{f3}) 
{\it 
Assume that $\# \Lambda_X = 1$ 
with $\Lambda_X = \{ \lambda_1 \}$. 
Then $\lambda_1 \in \mathbb{R}$. 
Moreover: 
}

\smallskip

(\roman{f3}-1) 
{\it 
$v_X = 1$ 
iff 
$\lambda_1 E \in {\cal O}_{G(K')^{\circ}}(X)$; 
}

\smallskip

(\roman{f3}-2) 
{\it 
$v_X = 2$ 
iff 
$\lambda_1 E + M_{1'} 
\in {\cal O}_{G(K')^{\circ}}(X)$; and  
}

\smallskip

(\roman{f3}-3) 
{\it 
$v_X = 3$ 
iff 
$\lambda_1 E + M_{2' 3} 
\in {\cal O}_{G(K')^{\circ}}(X)$. 
}

\medskip

\noindent
(2) 
{\it 
For $X, Y \in \mathcal{J}_3(K')$, 
${\cal O}_{G(K^{\mathbb{C}})^{\circ}}(X) 
= {\cal O}_{G(K^{\mathbb{C}})^{\circ}}(Y)$ 
iff $\Lambda_X = \Lambda_Y$ and $v_X = v_Y$. 
For any $X \in \mathcal{J}_3(K')$, 
${\cal O}(X) := {\cal O}_{G(K')^{\circ}}(X) 
= {\cal O}_{G(K')}(X)$ and 
${\cal O}(X) \cap \mathcal{J}_3(\bC') \ne \emptyset$. 
}

\bigskip

By Proposition 0.1 (1), 
$\Lambda_X$ and $v_X$ 
are invariants on ${\cal O}_{G(\tilde{K})}(X)$, 
so that the Theorems 0.2 (2) and 0.3 (2) 
follow from Theorems 0.2 (1) and 0.3 (1), 
respectively. 
Hence, 
this paper is concentrated in proving 
Theorems 0.2 (1) and 0.3 (1) 
with Proposition 0.1.

Note that the second equality of Propositoin 0.1 (1) 
was obtained by N.~Jacobson \cite[Lemma 1]{Jn1959}   
in a more general setting 
(cf. \cite[p.159, Proposition 5.9.4, \S 5.10]{SV}). 
In \S 1, by Lemma 1.2, 
it appears that 
the characteristic polynomial 
$\Phi_X(\lambda)$ of $X$ 
equals the generic minimal polynomial of $X$ 
defined by N.~Jacobson \cite[p.358 (5)]{Jn1968}. 
By Lemma 1.6 (3), 
it appears that $v_X$ 
equals the degree of 
N.~Jacobson \cite[p.389, Theorem 10]{Jn1968}'s 
minimal polynomial for 
$X \in \mathcal{J}_3(\tilde{K})$ 
with respect to the Jordan product.

\bigskip

\begin{center}
{\bf Contents}
\end{center}

\medskip

1. Preliminaries and Proposition 0.1 (1) and (2). 

2. Proposition 0.1 (3) and (4) (\roman{f1}). 

3. Theorems 0.2 and 0.3 in (1) (\roman{f1}, \roman{f2}). 

4. Proposition 0.1 (4) (\roman{f2}) 
and Theorems 0.2 and 0.3 in (1) (\roman{f3}).

Acknowledgement.

Reference.

\bigskip

\begin{center}
{\bf 1. Preliminaries 
and Proposition 0.1 (1) and (2).} 
\end{center} 

Let 
$i, i+1, i+2 \in \{1, 2, 3 \}$ 
be the indices counted modulo $3$. 
Then 

\begin{eqnarray*} 
& &
\left\{
\begin{array}{lll}
E_i \circ E_i = E_i, 
& E_i \circ E_{i+1} = 0, \\
E_i \circ F_i(x) = 0, 
& E_i \circ F_j(x) = \frac{1}{2} F_j(x)~(i \neq j),\\
F_i(x) \circ F_i(y) 
= (x | y) (E_{i+1} + E_{i+2}),
& F_i(x) \circ F_{i+1}(y) 
= \frac{1}{2} F_{i+2}(\overline{xy}); 
\end{array} \right. \\ 
& &
\left\{
\begin{array}{lll}
E_i \times E_i=0, 
& E_i \times E_{i+1} = \frac{1}{2}E_{i+2}, \\
E_i \times F_i(x) = - \frac{1}{2}F_i(x), 
& E_i \times F_j(x) = 0~(i \neq j), \\
F_i(x) \times F_i(y) 
= -(x | y)E_i,
& F_i(x) \times F_{i+1}(y) 
= \frac{1}{2}F_{i+2}(\overline{xy})  
\end{array}\right. \\
\end{eqnarray*}

\noindent 
for any $x, y \in \tilde{K}$. 
And 
 
\begin{eqnarray*}
& &M_1(x) \times M_1(y) 
= \sqrt{-1} \{ 
(x | 1) (y | 1) - (x | y) \} E_1, \\ 
& &M_1(x) \times M_{2 3}(y) 
= - \frac{1}{2} 
\{ 
F_2(\sqrt{-1}(\overline{x} - (x | 1))\overline{y}) 
+ F_3 (y (\overline{x} - (x | 1))) \}, \\
& &M_{23}(x) \times M_{23}(y) = (x|y) M_1,~
M_1 = M_{2 3}^{\times 2};~{\rm and} \\ 
& &M_{1'}(x) \times M_{1'}(y) 
= \{ (x|y) - (x|1)(y|1) \} E_1, \\
& &M_{1'}(x) \times M_{2' 3}(y) 
= \frac{1}{2} 
\{ F_2(- (\overline{x} \sqrt{-1} e_4) \overline{y} 
+ (x | 1) \sqrt{-1} e_4 \overline{y}) \\
& &~~~~~~~~~~~~~~~~~~~~~~~~~~~~
+ F_3(-(y \sqrt{-1}e_4)(\overline{x} \sqrt{-1}e_4) 
+ (x|1) y) \}, \\ 
& &M_{2 3'}(x) \times M_{2 3'}(y) 
= (x|y) M_{1'},~~M_{1'} = M_{2' 3}^{\times 2}. 
\end{eqnarray*}

\noindent

Let denote 
$\mathbb{X}(r; x) 
:= \sum_{i=1}^3 r_i E_i 
+ \sum_{i=1}^3 F_i(x_i)$ 
for any 
$r = (r_1, r_2, r_3) \in \mathbb{F}_3$ 
and  
$x = (x_1, x_2, x_3) \in \tilde{K}_3$. 
If $Y = \mathbb{X}(r; x) \in \mathcal{J}_3(\tilde{K})$, 
put $(Y)_{E_i}:= (Y | E_i) = r_i$ 
and $(Y)_{F_i}:= (Y | F_i(1))/2 = x_i$.

\bigskip

{\sc Lemma 1.1}. (1) 
{\it 
Let $i, i+1, i+2 \in \{1, 2, 3 \}$ 
be counted modulo $3$. 
Then 
}

\begin{eqnarray*} 
& &\mathbb{X}(r; x) \times \mathbb{X}(s; y) 
= \frac{1}{2} \sum_{i=1}^3 
\{ 
(r_{i+1} s_{i+2} + s_{i+1} r_{i+2} 
- 2 (x_i | y_i)) E_i \\
& &~~~~~~~~~~~~~~~~~~~~~~~~~~~~~~~~~~~
+ F_i(\overline{x_{i+1} y_{i+2} + y_{i+1} x_{i+2}} 
- r_i y_i - s_i x_i) \}; \\
& &
(\mathbb{X}(r; x) | \mathbb{X}(s; y)) 
= \sum_{i=1}^3 (r_i s_i + 2(x_i | y_i)); \\
& &
(\mathbb{X}(r; x) | \mathbb{X}(s; y) | \mathbb{X}(u; z)) 
= (\mathbb{X}(r; x) 
| \mathbb{X}(s; y) \times \mathbb{X}(u; z)) \\
& &~~~~~~~= 
\sum_{i=1}^3 
\{ 
\frac{r_i}{2} ( s_{i+1} u_{i+2} + u_{i+1} s_{i+2} ) 
+ (\overline{x_i} | y_{i+1} z_{i+2} + z_{i+1} y_{i+2}) \\
& &~~~~~~~~~~~~~~~~~~~~~~~~~~~~~~~~~~~~~~~~~~~~~~~
- r_i (y_i | z_i) - s_i (z_i | x_i) - u_i (y_i | x_i) \}; \\
& &\mathrm{det}(\mathbb{X}(r; x)) 
= r_1 r_2 r_3 + 2 (\overline{x_i} | x_{i+1} x_{i+2}) 
- \sum_{j = 1}^3 r_j N(x_j)~{\rm for}~i \in \{ 1, 2, 3 \}.   
\end{eqnarray*}

(2) 
{\it 
For $X, Y, Z \in \mathcal{J}_3(\tilde{K})$, 
all of $X \circ Y$, $(X | Y)$, $X \times Y$, 
$(X \circ Y | Z)$ and $(X | Y | Z)$ 
are symmetric. 
And $(X | Y)$ is non-degenerate. 
}

\smallskip

(3) 
{\it 
$2 E \times X = \mathrm{tr}(X) E - X  
= \varphi_{X} (\mathrm{tr}(X))$. 
Especially,  
$E^{\times 2} = E$ 
and 
$2 E \times X^{\times 2} 
= \mathrm{tr}(X^{\times 2}) E - X^{\times 2} 
= \varphi_{X^{\times 2}}( \mathrm{tr}(X^{\times 2}) )$.
}

\smallskip

(4) $(X | Y | E) = \mathrm{tr}(X \times Y) 
= \frac{1}{2}(\mathrm{tr}(X) \mathrm{tr}(Y) 
- (X | Y))$.

\medskip

{\it Proof.} 
(1) follows from the definitions except 
the 3rd equality, 
which is proved by \cite[p.15, 3.5 (7)]{Dd1978} 
as follows:  

\begin{eqnarray*}
& &(\mathbb{X}(r; x) | \mathbb{X}(s; y) | \mathbb{X}(u; z)) 
= \sum_{i=1}^3 \{ 
u_i ( r_{i+1} s_{i+2} + s_{i+1} r_{i+2} ) / 2 \\
&+& (\overline{z_i} | x_{i+1} y_{i+2} + y_{i+1} x_{i+2}) 
- u_i (x_i | y_i) - r_i (y_i | z_i) - s_i (x_i | z_i) \} \\
&=& \sum_{i=1}^3 \{ 
(u_{i+2} r_{i+3} s_{i+4} + u_{i+1} s_{i+2} r_{i+3}) / 2 \\ 
&+& (\overline{x_{i+3}}~\overline{z_{i+2}} | y_{i+4}) 
+ (\overline{z_{i+1}}~\overline{x_{i+3}} | y_{i+2}) 
- u_i (x_i | y_i) - r_i (y_i | z_i) - s_i (x_i | z_i) \} \\
&=& \sum_{i=1}^3 \{ 
r_i ( s_{i+1} u_{i+2} + u_{i+1} s_{i+2} ) / 2 
+ (\overline{x_i} | y_{i+1} z_{i+2} + z_{i+1} y_{i+2}) \\
&-& r_i (y_i | z_i) - s_i (z_i | x_i) - u_i (y_i | x_i) \} \\
&=& (\mathbb{X}(r; x) | \mathbb{X}(s; y) \times \mathbb{X}(u; z)). 
\end{eqnarray*}

\noindent 
(2) follows from the definitions or (1). 
(3) follows from direct computations. 
(4) follows from the definitions of $(X|Y|Z)$ and $X \times Y$.~\qed

\bigskip

For $X \in \mathcal{J}_3(\tilde{K})$, 
put 
$\Delta_X (\lambda) 
:= - \frac{1}{2} \{ 
3 \lambda^2 - 2 \mathrm{tr}(X) \lambda 
+ \mathrm{tr}(X)^2 - 2(X | X) \}$, 
which values in $\mathbb{F}$ (or ${\mathbb{C}}$) 
if $\lambda \in \mathbb{F}$ (resp. ${\mathbb{C}}$). 

\bigskip

{\sc Lemma 1.2.} 
(1) 
{\it 
$\Phi _X(\lambda) 
= \lambda^3 - \mathrm{tr}(X) \lambda^2 
+ \mathrm{tr}(X^{\times 2})\lambda - \mathrm{det}(X)$ 
with 
$\mathbb{F} \ni 
\mathrm{tr}(X) = \lambda_1 + \lambda_2 + \lambda_3$, 
$\mathrm{tr}(X^{\times 2}) 
= \lambda_1 \lambda_2 
+ \lambda_2 \lambda_3 + \lambda_3 \lambda_2$, 
$\mathrm{det}(X) = \lambda_1 \lambda_2 \lambda_3$ 
if $\Lambda_X 
= \{ \lambda_1, \lambda_2, \lambda_3 \} 
\subset {\mathbb{C}}$ 
for $X \in \mathcal{J}_3(\tilde{K})$. 
}

\smallskip

(2) 
{\it 
$\Phi_X'(\lambda) 
= 3 \lambda^2 - 2 \mathrm{tr}(X) \lambda 
+ \mathrm{tr}(X^{\times 2}) 
= \mathrm{tr}(\varphi_X(\lambda)^{\times 2})
= - 2 \Delta_X(\lambda) 
- \frac{1}{2} \{ \mathrm{tr}(X)^2 
- 3 (X | X) \}$. 
}

\smallskip

(3) 
{\it 
Put 
$\mathcal{M}(\tilde{K}) 
:= \{ X \in \mathcal{J}_3(\tilde{K})_0 
|~X \ne 0,~\Phi_X(\lambda) = \lambda^3 \}$. 
Then 
$\mathcal{M}(\tilde{K}) 
= \{ X \in \mathcal{J}_3(\tilde{K})_0 
|~X \ne 0, 
\mathrm{tr}(X^{ \times 2}) = \mathrm{det}(X) = 0 \} 
= \mathcal{M}_1(\tilde{K}) \cup \mathcal{M}_{2 3}(\tilde{K})$ 
with 
$\mathcal{M}_1(\tilde{K}) \cap 
\mathcal{M}_{2 3}(\tilde{K}) = \emptyset$.  
And 
$\{ X \in \mathcal{J}_3(\tilde{K}) 
|~\# \Lambda_X = 1 \} = 
\mathbb{F} E \oplus 
(\{ 0 \} \cup \mathcal{M} (\tilde{K}))$.  
}

\medskip

{\it Proof.} 
(1) 
$\Phi_X(\lambda) = \frac{1}{3}
( \lambda E - X | \lambda E - X | \lambda E - X)$, 
which equals the required one by 
Lemma 1.1 (2, 3) and 
$\Phi_X (\lambda) \equiv 
(\lambda - \lambda_1) 
(\lambda - \lambda_2) (\lambda -\lambda_3)$. 

(2) 
The first equality folllows from (1). 
By Lemma 1.1 (3), 
$\varphi_X(\lambda)^{\times 2} 
= (\lambda E - X)^{\times 2} 
= \lambda^2 E - (\mathrm{tr}(X) E - X) \lambda + X^{\times 2}$, 
so that 
$\mathrm{tr} (\varphi_X(\lambda)^{\times 2}) 
= 3 \lambda^2 - 2 \mathrm{tr}(X) \lambda + \mathrm{tr}(X^{\times 2})$. 
And 
$3 \lambda^2 - 2 \mathrm{tr}(X) \lambda + \mathrm{tr}(X^{\times 2}) 
= - 2 \Delta_X(\lambda) - \frac{1}{2} \{ \mathrm{tr}(X)^2 - 3 (X | X) \}$ 
by the second equality of Lemma 1.1 (4).

(3) 
The first claim follows from (1). 
For $X \in \mathcal{J}_3(\tilde{K})$, 
put 
$X_0 := X - \frac{1}{3} \mathrm{tr}(X) E 
\in \mathcal{J}_3(\tilde{K})_0$. 
Then 
$X = \frac{1}{3} \mathrm{tr}(X) E + X_0$, 
so that 
$\mathcal{J}_3(\tilde{K}) 
= \mathbb{F} E \oplus \mathcal{J}_3(\tilde{K})_0$. 
If $\Phi_X(\lambda) 
= \Pi_{i = 1}^3 (\lambda - \lambda_i)$, 
then 
$\Phi_{X_0}(\lambda) 
= \mathrm{det}((\lambda + \frac{1}{3} \mathrm{tr}(X))E - X) 
= \Pi_{i=1}^3 (\lambda + \frac{1}{3} \mathrm{tr}(X) - \lambda_i)$, 
so that 
$\Phi_{X_0}(\lambda) = \lambda^3 
\Leftrightarrow 
\frac{1}{3} \mathrm{tr}(X) - \lambda_i = 0 
(i = 1, 2, 3) 
\Leftrightarrow 
\lambda_1 = \lambda_2 = \lambda_3 
\Leftrightarrow \# \Lambda_X = 1$, 
because of $\mathrm{tr}(X) = \sum_{i =1}^3 \lambda_i$ 
by (1). 
Hence, 
$\{ X \in \mathcal{J}_3(\tilde{K}) 
|~\# \Lambda_X = 1 \} = 
\mathbb{F} E \oplus 
(\{ 0 \} \cup \mathcal{M} (\tilde{K}))$.~\qed 

\bigskip

Let $V$ be an 
$\mathbb{F}$-algebra 
with the multiplication $x y$ 
of $x, y \in V$. 
For $x \in V$, put an 
$\mathbb{F}$-linear endomorphism on $V$,  
$L_x: V \rightarrow V; y \mapsto x y$, 
as the left translation by $x$. 
And put the automorphism group of $V$ as follows: 

\[
\mathrm{Aut}(V) := 
\{ \alpha \in GL_{\mathbb{F}}(V) 
|~\alpha(x y) = (\alpha x) (\alpha y);
~x, y \in V \}. 
\]

\bigskip

{\sc Lemma 1.3.} 
(1) 
{\it 
Let $V$ be an $\mathbb{F}$-algbra. 
Assume that 
$\alpha \in \mathrm{Aut}(V)$. 
Then 
$\mathrm{trace} (L_{(\alpha x)}) 
= \mathrm{trace} (L_x)$,
~$\mathrm{det} (L_{(\alpha x)}) 
= \mathrm{det} (L_x)$ for all $x \in V$. 
If moreover $V$ 
admits the identity element $e$, 
then $\alpha e = e$. 
}

\smallskip 

(2) 
{\it 
Let $L^{\circ}_X$ 
and $L^{\times}_X$ 
be the left translations by 
$X \in \mathcal{J}_3(\tilde{K})$ 
on $\mathcal{J}_3(\tilde{K})$ 
with respect to the product $\circ$ 
and the cross product $\times$, 
respectively. 
Then 
$\mathrm{trace}(L^{\circ}_X) 
= (d_K + 1)~\mathrm{tr}(X)$ 
and 
$\mathrm{trace}(L^{\times}_X) 
= \frac{- 1}{2} d_K~\mathrm{tr}(X)$. 
}

\medskip

{\it Proof.} 
(1) 
For $x, y \in V$ and $\alpha \in G(V)$, 
$L_{(\alpha x)} y = (\alpha x) y 
= \alpha (x (\alpha^{-1} y)) 
= (\alpha L_x \alpha^{-1}) y$, 
i.e. 
$L_{(\alpha x)} = \alpha L_x \alpha^{-1}$, 
so that 
$\mathrm{trace} (L_{(\alpha x)}) 
= \mathrm{trace} (L_x)$ 
and $\mathrm{det}(L_{(\alpha x)}) 
= \mathrm{det} (L_x)$ 
as an $\mathbb{F}$-linear endomorphism on $V$. 
Assume that $e x = x e = x$ for any $x \in V$. 
Take $\alpha \in \mathrm{Aut}(V)$. 
Then 
$(\alpha e)(\alpha x) 
= (\alpha x)(\alpha e) = \alpha x$, 
so that 
$(\alpha e) y 
= y (\alpha e) = y$ for all $y \in V$. 
In particular, 
$\alpha e = (\alpha e) e = e$. 

(2) 
$\{ E_i, F_i(e_j/\sqrt{2}) 
|~i = 1, 2, 3; j = 0, \cdots, d_K-1 \}$ 
forms an orthonormal basis of 
$(\mathcal{J}_3(K^{{\mathbb{C}}}), (*|*))$ 
by Lemma 1.1 (1). 
And $L^{\circ}_X$ and $L^{\times}_X$ 
can be identified with 
a ${\mathbb{C}}$-linear endomorphism on 
$\mathcal{J}_3(K^{{\mathbb{C}}}) 
= \mathcal{J}_3(\tilde{K})$ or 
${\mathbb{C}} \otimes \mathcal{J}_3(\tilde{K})$. 
By Lemma 1.1 (1, 2), 
$\mathrm{trace}(L^{\circ}_X) 
= \sum_{i=1}^3 \{ (X \circ E_i | E_i) 
+ \frac{1}{2} \sum_{j=0}^{d_K-1} 
(X \circ F_i(e_j) | F_i(e_j)) \} 
= \sum_{i=1}^3 \{ 
(X | E_i \circ E_i) 
+ \frac{1}{2} \sum_{j=0}^{d_K-1} 
(X | F_i(e_j) \circ F_i(e_j)) \} 
= \sum_{i=1}^3 \\
\{ (X | E_i) + \frac{1}{2} 
\sum_{j=0}^{d_K-1} (X | E_{i+1} + E_{i+2}) \}  
= (d_K +1)~\mathrm{tr}(X)$; 
and $\mathrm{trace}(L_X^{\times}) 
= \sum_{i=1}^3 \{ (X \times E_i, E_i) 
+ \frac{1}{2} \sum_{j=0}^{d_K-1} 
(X \times F_i(e_j) | F_i(e_j)) \} 
= \sum_{i=1}^3 \{ (X | E_i \times E_i) 
+ \frac{1}{2} \sum_{j=0}^{d_K-1} 
(X | F_i(e_j) \times F_i(e_j)) \} 
= \sum_{i=1}^3 \frac{1}{2} 
\sum_{j=0}^{d_K-1} (X | -E_i) 
= \frac{- 1}{2} d_K \mathrm{tr}(X)$. 
\qed

\bigskip

{\it Proof of Proposition 0.1 (1).} 
The first claim follows from Lemma 1.3 (1)(2). 
For the second claim, 
since $\mathrm{det}(X)$ 
is defined by $X \circ X$, 
$\mathrm{tr}(X)$ and $E$, 
the first equality is recognized as 
the inclusion $\subseteqq$. 
By $\Phi_X(\lambda) = \mathrm{det}(\lambda E - X)$, 
the 2nd equality is recognized as 
the inclusion $\subseteqq$.  
By Lemmas 1.1 (4) and 1.2 (1), 
the 3rd equality is recognized 
as the inclusion $\subseteqq$. 
By polarizing $3 {\rm det}(X) = (X|X|X)$ 
with Lemma 1.1 (2), 
the 4th equality is recognized as the inclusion 
$\subseteqq$. 
Assume that 
$\alpha \in \mathrm{GL}_{\mathbb{F}}
(\mathcal{J}_3(\tilde{K}))$ 
and 
$(\alpha X) \times (\alpha Y) = \alpha (X \times Y)$ 
for all 
$X, Y \in \mathcal{J}_3(\tilde{K})$. 
By Lemma 1.3, 
$\mathrm{tr}(\alpha X) = \mathrm{tr}(X)$.  
By Lemma 1.1 (4), 
$(X|Y) = \mathrm{tr}(X) \mathrm{tr}(Y)
- 2 \mathrm{tr}(X \times Y)$, 
so that 
$(\alpha X | \alpha Y) = (X | Y)$. 
By the definition of $\times$, 
$(X \circ Y | Z)  
= (X \times Y | Z) 
+ ( \mathrm{tr}(X) (Y | Z) 
+ \mathrm{tr}(Y) (X | Z) 
- (\mathrm{tr}(X) \mathrm{tr}(Y) 
- (X | Y)) \mathrm{tr}(Z))/2$, 
so that 
$((\alpha X) \circ (\alpha Y) | \alpha Z) 
= (X \circ Y | Z)$ 
for all $X, Y, Z \in \mathcal{J}_3(\tilde{K})$. 
By Lemma 1.1 (2), 
$\alpha^{-1}((\alpha X) \circ (\alpha Y)) 
= X \circ Y$, that is, $\alpha \in G(\tilde{K})$. 
Hence, all of the equations of the second claim follow. 
The last claim follows from these equations. 

\medskip

{\it Proof of Proposition 0.1 (2).} 
Note that 
$\tau \gamma = \gamma \tau$, 
$\gamma E_i = E_i$, 
$\gamma E = E$ 
and 
$\mathrm{det}( \gamma \mathbb{X}(r; x) ) 
= \mathrm{det}( \mathbb{X}(r; x) )$ 
by Lemma 1.1 (1), 
so that 
$\gamma \in G(\tilde{K})^{\tau}_{E_1, E_2, E_3}$. 
By Proposition 0.1 (1), 
the last claim follows. 
For the first claim, put 
$< X | Y > := (\tau X | Y) \in \mathbb{F}$ 
for $X, Y \in \mathcal{J}_3(\tilde{K})$, 
which defines a positive-definite 
symmetric (or hermitian) 2-form on 
$\mathcal{J}_3(K')$ 
(resp. $\mathcal{J}_3(K^{\mathbb{C}})$) 
over $\mathbb{R}$ (resp. $\mathbb{C}$) 
by Lemma 1.1 (1). 
For $\alpha \in G(\tilde{K})$, 
$\alpha^* \in GL_F(\mathcal{J}_3(\tilde{K}))$ 
is defined such that 
$<\alpha X | Y> = < X | \alpha^* Y>$ 
for all 
$X, Y \in \mathcal{J}_3(\tilde{K})$. 
By (1), 
$< X | \alpha^* Y > 
= (\tau \alpha X | Y) 
= \tau (\alpha X | \tau Y) 
= \tau (X | \alpha^{-1} \tau Y) 
= <X | \tau \alpha^{-1} \tau Y>$, 
so that 
$\alpha^* = \tau \alpha^{-1} \tau 
\in G(\tilde{K})$ 
because of (1) by 
$\mathrm{det}(\alpha^* X) 
= \tau \mathrm{det}(\alpha^{-1} \tau X) 
= \tau^2 \mathrm{det}(X) 
= \mathrm{det}(X)$ 
and $\alpha^* E = \tau \alpha^{-1} \tau E = E$. 
Then 
$G(\tilde{K}) \cong G(\tilde{K})^{\tau} \times 
\bR^{}$ as a polar decomposition 
of C.~Chevalley \cite[p.201]{Cc} 
(resp. \cite[p.450, Lemma 2.3]{Hs}), 
so that $G(\tilde{K})^{\tau}$ 
is a maximal compact subgroup of $G(\tilde{K})$.~\qed 

\bigskip

For $i \in \{ 1, 2, 3 \}$ and $a \in \tilde{K}$, 
put $B_i(a): 
\mathcal{J}_3(\tilde{K}) 
\longrightarrow 
\mathcal{J}_3(\tilde{K}); 
\mathbb{X}(r; x) \mapsto 
\mathbb{X}(s; y)$ 
such that 
$s_i := 0, s_{i+1} := 2(a | x_i), 
s_{i+2} := - 2(a | x_i), 
y_i := - (r_{i+1} - r_{i+2}) a, 
y_{i+1} := - \overline{x_{i+2} a},
y_{i+2} := \overline{a x_{i+1}}$, 
where 
$i, i+1, i+2 \in \{ 1, 2, 3 \}$ 
are counted modulo $3$. 
Then 
$\mathrm{exp}(t B_i(a)) 
\in (G(\tilde{K})_{E_i})^{\circ}$ 
for $t \in \mathbb{F}$. 
In fact, 
$s_i = s_{i+1} = s_{i+2} 
= y_i = y_{i+1} = y_{i+2} = 0$ 
if $\mathbb{X}(r; x) = E_i$ or $E$. 
Put $X = \mathbb{X}(r; x)$. 
By Lemma 1.1 (1), 
$(B_i(a) X | X | X)  
= (B_i(a) X | X^{\times 2}) 
= 2 \{ 
(a | x_i) (r_{i+2} r_i - N(x_{i+1}) 
- r_i r_{i+1} + N(x_{i+2})) 
- (r_{i+1} - r_{i+2}) 
(a | \overline{x_{i+1} x_{i+2}} - r_i x_i) 
- (\overline{x_{i+2} a} 
| \overline{x_{i+2} x_i} - r_{i+1} x_{i+1}) 
+ (\overline{a x_{i+1}} 
| \overline{x_i x_{i+1}} - r_{i+2} x_{i+2}) 
\} = 0$, 
so that 
$\mathrm{exp}(t B_i(a)) 
\in (G(\tilde{K})_{E_i})^{\circ}$ 
for all $t \in \mathbb{F}$ 
by Proposition 0.1 (1), as required.  
Note that $B_i(a)$ 
is nothing but $\tilde{A}_i^a$ given 
in H.~Freudenthal~\cite[(5.1.1)]{Fh1951}.

For $\nu \in \{ 1, \sqrt{-1} \}$, 
put 
$C_{\nu}(t) := (e^{\nu t} + e^{- \nu t})/2$, 
$S_{\nu}(t) := (e^{\nu t} - e^{- \nu t})/(2 \nu)$ 
as $\mathbb{F}$-valued functions of 
$t \in \mathbb{F}$. 
Then 
$(C_{\nu}(t), S_{\nu}(t)) =  (\cosh(t), \sinh (t))$ 
or $(\cos(t), \sin (t))$ if $\nu = 1$ or $\sqrt{-1}$, 
respectively. 
Note that 

\begin{eqnarray*}
& &\tau C_{\nu}(t) = C_{\nu}(\tau t),
~\tau S_{\nu}(t) = S_{\nu}(\tau t), \\
& &C_{\nu}(t_1) C_{\nu}(t_2) 
+ \nu^2 S_{\nu}(t_1) S_{\nu}(t_2) 
= C_{\nu}(t_1 + t_2), \\
& &C_{\nu}'(t) = \nu^2 S_{\nu}(t),
~S_{\nu}'(t) = C_{\nu}(t), \\
&&S_{\nu}(t_1) C_{\nu}(t_2) + C_{\nu}(t_1) S_{\nu}(t_2) 
= S_{\nu}(t_1 + t_2), \\
& &C_{\nu}'(0) = 0,
~S_{\nu}'(0) = 1,
~C_{\nu}(2t) = 1 + 2 \nu^2 S^2_{\nu}(t). 
\end{eqnarray*}

For $i \in \{ 1, 2, 3 \}$, $t \in \mathbb{F}$, 
$a \in \tilde{K}$ and $\nu \in \{ 1, \sqrt{-1} \}$, 
put 
$\beta_i(t; a, \nu): 
\mathcal{J}_3(\tilde{K}) 
\longrightarrow 
\mathcal{J}_3(\tilde{K}); 
\mathbb{X}(r; x) 
\mapsto \mathbb{X}(s; y)$ 
such that

\[ 
\left\{
\begin{array}{ccl}
s_i &:=& r_i,\\
s_{i+1} &:=& 
\frac{r_{i+1} + r_{i+2}}{2} 
+ \frac{r_{i+1} - r_{i+2}}{2} C_{\nu}(2t) 
+ (a | x_i) S_{\nu}(2t),\\
s_{i+2} &:=& 
\frac{r_{i+1} + r_{i+2}}{2} 
- \frac{r_{i+1} - r_{i+2}}{2} C_{\nu}(2t) 
- (a | x_i) S_{\nu}(2t), \\
y_i &:=& 
x_i 
- a \frac{r_{i+1} - r_{i+2}}{2} S_{\nu}(2t) 
- 2 a (a | x_i) S_{\nu}^2(t),\\
y_{i+1} &:=& 
x_{i+1} C_{\nu}(t) - \overline{x_{i+2} a} S_{\nu}(t),\\
y_{i+2} &:=& 
x_{i+2} C_{\nu}(t) + \overline{ax_{i+1}} S_{\nu}(t). 
\end{array} \right.
\]

\noindent 
For $c \in \mathbb{F}$, 
put 
$\mathcal{S}_1(c, \tilde{K}) 
:= \{ x \in \tilde{K} 
|~N(x) = c \}$, 
which is said to be 
{\it a generalized sphere} 
\cite[p.42, (3.7)]{Hfr1990} 
{\it of first kind over} 
$\mathbb{F}$.

\bigskip

{\sc Lemma 1.4.} 
(1) 
(\roman{f1}) 
{\it 
Assume that 
$i \in \{ 1, 2, 3 \}$, 
$\nu \in \{ 1, \sqrt{-1} \}$ 
and 
$a \in \mathcal{S}_1(- \nu^2, \tilde{K})$. 
Then 
$\beta_i(t; a, \nu) 
= \mathrm{exp}(t B_i(a)) 
\in (G(\tilde{K})_{E_i})^{\circ}$ 
for $t \in \mathbb{F}$ such that 
$\beta_i(t; a, \nu) \tau 
= \tau \beta_i(\tau t; \tau a, \nu)$ 
for all $t \in \mathbb{F}$. 
Especially, 
$\sigma_i 
:= \beta_i(\pi; 1, \sqrt{-1}) 
\in ((G(\tilde{K})_{E_i}^{\tau})^{\circ})
_{E_{i+1}, E_{i+2}}$. 
}

\smallskip

(\roman{f2}) 
{\it 
For $i \in \{ 1, 2, 3 \}$, put 
$\hat{\beta}_i 
:= \beta_i(\frac{\pi}{2}; 1, \sqrt{-1})$. 
Then 
$\hat{\beta}_i 
\in (G(\tilde{K})_{E_i}^{\tau})^{\circ}$ 
such that 
$\hat{\beta}_i X 
= r_i E_i + r_{i+2} E_{i+1} 
+ r_{i+1} E_{i+2}  
+ F_{i} (- \overline{x_i}) 
+ F_{i+1} (- \overline{x_{i+2}}) 
+ F_{i+2} (\overline{x_{i+1}})$ 
if $X = \mathbb{X}(r; x) \in \mathcal{J}_3(\tilde{K})$.  
Especially, for any permutation 
$\mu = (\mu_1, \mu_2, \mu_3)$ 
of the triplet $(1, 2, 3)$, 
there exists 
$\hat{\beta} \in (G(\tilde{K})^{\tau})^{\circ}$ 
such that 
$\hat{\beta} ( \sum_{j = 1}^3 r_j E_j ) = 
\sum_{j = 1}^3 r_{\mu_j} E_{\mu_j}$ 
for all $r_i \in \mathbb{F}~(i = 1, 2, 3)$.
}

\smallskip

(\roman{f3}) 
{\it 
Put 
$B_{2 3} := B_2(\sqrt{-1}) - B_3(1)$,  
$B_{2 3'} := B_2(1) - B_3(\sqrt{-1} e_4)$, 
$B_{2' 3} := B_2(-\sqrt{-1}e_4) - B_3(1)$, 
$\beta_{2 3}(t) := \mathrm{exp}(t B_{2 3})$, 
$\beta_{2 3'}(t) := \mathrm{exp}(t B_{2 3'})$, 
$\beta_{2' 3}(t) := \mathrm{exp}(t B_{2' 3})$. 
Then 
$\beta_{2 3}(t) \in (G(K^{\mathbb{C}})^{\circ})_{M_1}$ 
and 
$\beta_{2 3}(t) M_{2 3}(x) 
= 2 t (x | 1) M_1 + M_{2 3}(x)$ 
$(x \in K^{\mathbb{C}}$, $t \in \mathbb{C})$. 
And 
$\beta_{2 3'}(t), \beta_{2' 3}(t) 
\in (G(K')^{\circ})_{M_{1'}}$ 
such that 
$\beta_{2 3'}(t) M_{2' 3}(x) 
= 2 t (\sqrt{-1} e_4 | x) M_{1'} + M_{2' 3}(x)$, 
$\beta_{2' 3}(t) M_{2' 3}(x) 
= 2 t (1 | x) M_{1'} + M_{2' 3}(x)$ 
$(x \in K'$, $t \in \mathbb{R})$. 
} 

\smallskip

(2) 
(\roman{f1}) 
{\it 
Let 
$\mathcal{S}_1(1, \tilde{K})^{\circ}$ 
be the connected component of 
$\mathcal{S}_1(1, \tilde{K})$ 
containing $1 = e_0$ 
in $\tilde{K}$. 
And 
$O(\tilde{K}) := 
\{ \alpha \in GL_F (\tilde{K}) 
|~N(\alpha x) = N(x) \}$. 
Then 
$\mathcal{S}_1(1, \tilde{K}) 
= {\cal O}_{O(\tilde{K})} (e_0) 
= \mathcal{S}_1(1, \tilde{K})^{\circ} 
\cup (- \mathcal{S}_1(1, \tilde{K})^{\circ})$. 
Especially, 
$\mathcal{S}_1(1, \tilde{K}) 
= \mathcal{S}_1(1, \tilde{K})^{\circ} 
= - \mathcal{S}_1(1, \tilde{K})^{\circ}$ 
when 
$\tilde{K} = \bH', \bO'; \bC^{\mathbb{C}}, 
\bH^{\mathbb{C}}, \bO^{\mathbb{C}}$.
}

\smallskip

(\roman{f2}) 
{\it 
For 
$a \in \mathcal{S}_1(1, \tilde{K})$ 
and $i \in \{ 1, 2, 3 \}$, 
put 
$\delta_i(a) \in 
{\rm End}_{\mathbb{F}}(\mathcal{J}_3(\tilde{K}))$ 
with 
$\mathbb{X}(s; y) 
:= \delta_i(a) \mathbb{X}(r; x)$ 
such that 
$s_i := r_i, 
s_{i+1} := r_{i+1}, 
s_{i+2} := r_{i+2}, 
y_{i} := a x_{i} a, 
y_{i+1} := \overline{a} x_{i+1}, 
y_{i+2} := x_{i+2} \overline{a}$. 
Then 
$\delta_i(a) \in 
((G(\tilde{K})_{E_i})^{\circ})
_{E_{i+1}, E_{i+2}}$ 
such that 
$\delta_i(a) \sigma_i 
= \sigma_i \delta_i(a) 
= \delta_i(-a)$ 
and 
$\delta_i(a) \tau 
= \tau \delta_i(\tau a)$. 
Especially, 
$\delta_i(a) \in 
(G(\tilde{K})_{E_1, E_2, E_3})^{\circ}$ 
when 
$\tilde{K} = \bH', \bO'; \bC^{\mathbb{C}}, 
\bH^{\mathbb{C}}, \bO^{\mathbb{C}}$.
}

\smallskip

(\roman{f3}) 
{\it 
Assume that 
$d_{\tilde{K}} \leqq 4$. 
For 
$a \in \mathcal{S}_1(1, \tilde{K})$ 
and $i \in \{ 1, 2, 3 \}$, 
put 
$\beta_i(a) \in 
{\rm End}_{\mathbb{F}}(\mathcal{J}_3(\tilde{K}))$ 
with 
$\mathbb{X}(s; y) := 
\beta_i(a) \mathbb{X}(r; x)$ 
such that 
$s_i := r_i, 
s_{i+1} := r_{i+1}, 
s_{i+2} := r_{i+2}, 
y_{i} := a x_{i} \overline{a}, 
y_{i+1} := a x_{i+1}, 
y_{i+2} := x_{i+2} \overline{a}$.  
Then 
$\beta_i(a) \in 
((G(\tilde{K})_{E_i})^{\circ})
_{E_{i+1}, E_{i+2}, F_i(1)}$ 
such that 
$\beta_i(a) \sigma_i 
= \sigma_i \beta_i(a) 
= \beta_i(-a)$. 
Especially, 
$\beta_i(a) \in 
(G(\tilde{K})_{E_1, E_2, E_3, F_1(1)})^{\circ}$ 
when 
$\tilde{K} = \bH'; 
\bC^{\mathbb{C}}, 
\bH^{\mathbb{C}}$.
}

\medskip

{\it Proof.} 
(1) (\roman{f1}) 
Put 
$\mathbb{X}(u; z) := 
\frac{d}{d t} 
\beta_i(t; a, \nu) 
\mathbb{X}(r; x) - B_i(a) \mathbb{X}(r; x)$. 
Then 
$u_i = z_i = 0$, 
$u_{i+1} 
= (\nu^2 + N(a)) ((r_{i+1}-r_{i+2}) S_{\nu}(2 t) 
+ 4 (a, x_i) S_{\nu}^2(t)) = 0 = - u_{i+2}$, 
$z_{i+1} 
= (\nu^2 + N(a)) x_{i+1} S_{\nu}(t) = 0$, 
$z_{i+2} = (\nu^2 + N(a)) x_{i+2} S_{\nu}(t) = 0$, 
i.e. 
$\frac{d}{d t} 
\beta_i(t; a, \nu) \mathbb{X}(r; x) 
= B_i(a) \mathbb{X}(r; x)$ ($t \in \mathbb{F}$) 
with 
$\beta_i(0; a, \nu) \mathbb{X}(r; x) = \mathbb{X}(r; x)$. 
Hence, 
$\beta_i(t; a, \nu) = \mathrm{exp}(t B_i(a))$, 
so that 
$\beta_i(t; a, \nu) \in (G(\tilde{K})_{E_i})^{\circ}$ 
and 
$\beta_i(t; a, \nu) \tau 
= \mathrm{exp}(t B_i(a)) \tau 
= \tau \mathrm{exp}((\tau t) B_i(\tau a)) 
= \tau \beta_i(t; a, \nu)$ 
for all $t \in \mathbb{F}$.  
Especially, 
$\sigma_i(t) := \beta_i(\pi t; 1, \sqrt{-1}) 
\in (G(\tilde{K})_{E_i}^{\tau})^{\circ}$ 
for all $t \in \mathbb{R}$ 
such that $\sigma_i = \sigma_i(1)$, 
$\sigma_i E_{i+1} = E_{i+1}$, 
$\sigma_i E_{i+2} = E_{i+2}$.

(\roman{f2}) 
The first claim follows from (\roman{f1}), 
so that the second claim follows.

(\roman{f3}) 
For $a, b \in K^{\mathbb{C}}$, 
$(B_2(a) - B_3(b)) M_1 
= M_{2 3}( - \sqrt{-1} \overline{a} - b)$. 
Then $B_{2 3} M_1 = 0$ 
by $a = \sqrt{-1}$ and $b = 1$. 
For $x \in K^{\mathbb{C}}$, 
$(B_2(a) - B_3(b)) M_{2 3}(x) 
= \mathbb{X}(-2((a | \sqrt{-1} \overline{x})+(b|x)), 
2(b|x), 2 (a | \sqrt{-1} \overline{x}); 
\overline{a x}+ \sqrt{-1}~\overline{b} x, 0, 0)$. 
In particular, 
$B_{2 3} M_{2 3}(x) = 2(1|x) M_1$. 
Hence, 
$\beta_{2 3}(t) \in (G(K^{\mathbb{C}})^{\circ})_{M_1}$ 
and 
$\beta_{2 3}(t) M_{2 3}(x) 
= 2 t (1 | x) M_1 + M_{2 3}(x)$ 
for $x \in K^{\mathbb{C}}$ and $t \in \mathbb{C}$. 
For $a, b \in K'$, 
$(B_2(a) - B_3(b)) M_{1'} 
= M_{2' 3}(\overline{a} \sqrt{-1} e_4 - b)$. 
Then 
$B_{2 3'} M_{1'} = B_{2' 3} M_{1'} = 0$ 
by $b = \overline{a} \sqrt{-1} e_4$ 
with 
$(a, b) = (1, \sqrt{-1} e_4), (-\sqrt{-1} e_4, 1)$. 
For $x \in K'$, 
$(B_2(a) - B_3(b)) M_{2' 3}(x) = 
\mathbb{X}(-2((a | -\sqrt{-1}e_4 \overline{x}) + (b|x)), 
2(b|x), 2(a|-\sqrt{-1}e_4 \overline{x}); 
\overline{a x} 
- \overline{(\sqrt{-1}e_4 \overline{x}) b}, 0, 0)$,  
so that 
$B_{2 3'} M_{2' 3}(x) 
= 2(\sqrt{-1}e_4|x) E_2 
- 2(1|\sqrt{-1}e_4 \overline{x}) E_3 
+ F_1(\overline{x} + \sqrt{-1} e_4 (x \sqrt{-1} e_4))$ 
and 
$B_{2' 3} M_{2' 3}(x) 
= 2 (1 | x) M_{1'}$. 
Put $x = p + q \sqrt{-1} e_4$ 
with $p, q \in \bH$, 
so that 
$\overline{x} = \overline{p} - q \sqrt{-1} e_4$. 
Then 
$(\sqrt{-1} e_4) (x \sqrt{-1} e_4) 
= \overline{p} + \overline{q} \sqrt{-1} e_4$, 
$\overline{x} - \sqrt{-1} e_4 (x \sqrt{-1} e_4) 
= - (q + \overline{q}) \sqrt{-1} e_4 
= 2 (\sqrt{-1} e_4 | x) \sqrt{-1} e_4$. 
And 
$B_{2 3'} M_{2' 3}(x) = 2(\sqrt{-1} e_4 | x) M_{1'}$. 
Hence, 
$\beta_{2 3'}(t), \beta_{2' 3}(t) 
\in (G(K^{\mathbb{C}})^{\circ})_{M_{1'}}$ 
such that 
$\beta_{2 3'}(t) M_{2' 3}(x) 
= 2 t (\sqrt{-1} e_4 |x) M_{1'} + M_{2' 3}(x)$ 
and 
$\beta_{2' 3}(t) M_{2' 3}(x) 
= 2 t (1 | x) M_{1'} + M_{2' 3}(x)$ 
for $x \in K'$ and $t \in \mathbb{R}$.

(2) 
(\roman{f1}) 
Since $\tilde{K}$ is a composition algebra, 
$L_a \in O(\tilde{K})$ 
for all $a \in \mathcal{S}_1(1, \tilde{K})$. 
Hence, 
$\mathcal{S}_1(1, \tilde{K}) 
= \{ L_a (e_0) |~a \in 
\mathcal{S}_1(1, \tilde{K}) \} 
= \mathcal{O}_{O(\tilde{K})}(e_0)$. 
Put 
$SO(\tilde{K}) := 
\{ \alpha \in O(\tilde{K}) 
|~{\rm det}(\alpha) = 1 \}$. 
When $d_{\tilde{K}} = 1$: 
$\mathcal{S}_1(1, \tilde{K}) = \{ \pm e_0 \}$, 
$\mathcal{S}_1(1, \tilde{K})^{\circ} 
= \{ e_0 \}$, 
$\mathcal{S}_1(1, \tilde{K}) 
= \mathcal{S}_1(1, \tilde{K})^{\circ} 
\cup (- \mathcal{S}_1(1, \tilde{K})^{\circ})$. 
When $d_{\tilde{K}} = 2, 4, 8$: 
$O(\tilde{K}) = 
SO(\tilde{K}) \cup SO(\tilde{K}) \epsilon$ 
with $\epsilon (e_0) = e_0$, 
$SO(\tilde{K}) = - SO(\tilde{K})$, 
so that 
$\mathcal{S}_1(1, \tilde{K}) 
= {\cal O}_{SO(\tilde{K})}(e_0) 
= - {\cal O}_{SO(\tilde{K})}(e_0) 
= - \mathcal{S}_1(1, \tilde{K})$. 
Since $SO(K^{\mathbb{C}})$ is connected, 

\[
\mathcal{S}_1(1, K^{\mathbb{C}})^{\circ} 
= \mathcal{S}_1(1, K^{\mathbb{C}}) 
= - \mathcal{S}_1(1, K^{\mathbb{C}})^{\circ}. 
\]

\noindent 
And 
$M_{d_{K'}}(\mathbb{R}) 
\supseteqq 
SO(K') \cong S(O(d_{K'}/2) \times O(d_{K'}/2)) 
\times \mathbb{R}^{(d_{K'}/2)^2}$. 
Put 

\[
1_n := \mathrm{diag}(1, \cdots, 1), 
1_n' := \mathrm{diag}(1_{n-1}, -1) \in M_n(K'). 
\]

\noindent
When $d_{K'}/2 = 2, 4$, 
$SO(K')$ 
admits just four connected components containing 
$\mathrm{id}_{K'}$, 
$\mathrm{diag}(1_{d_{K'}/2}', 1_{d_{K'}/2})$, 
$\mathrm{diag}(1_{d_{K'}/2}, 1_{d_{K'}/2}')$, 
$\mathrm{diag}(1_{d_{K'}/2}', 1_{d_{K'}/2}')$, 
so that 
$\mathcal{S}_1(1, K') 
= {\cal O}_{SO(K')}(e_0) 
= {\cal O}_{SO(K')^{\circ}}(e_0) 
= \mathcal{S}_1(1, K')^{\circ}$.

(\roman{f2}) 
For $a \in \mathcal{S}_1(1, \tilde{K})$, 
$\delta_i(a) \in 
GL_{\mathbb{F}}(\mathcal{J}_3(\tilde{K}))
_{E_1, E_2, E_3}$ 
with 
$\delta_i(a)^{-1}= \delta_i(\bar{a})$, 
$\delta_i(a) \tau = \tau \delta_i(\tau a)$ 
and $\delta_i(a) E = E$. 
By Lemma 1.1 (1), 
$\mathrm{det}(\delta_i(a) \mathbb{X}(r; x)) 
= r_1 r_2 r_3 
+ 2 ( \overline{a x_i a} | 
(\overline{a} x_{i+1}) (x_{i+2} \overline{a}) ) 
- r_i N(a x_i a) 
- r_{i+1} N(\overline{a} x_{i+1}) 
- r_{i+2} N(x_{i+2} \overline{a}) 
= r_1 r_2 r_3 
+ (2 (\bar{x}_i | x_{i+1} x_{i+2}) - r_i N(x_i)) N(a)^2 
- (r_{i+1} N(x_{i+1}) + r_{i+2} N(x_{i+2})) N(a)$ 
$= \mathrm{det}(\mathbb{X}(r; x))$. 
By Proposition 0.1 (1), 
$\delta_i(a) \in 
(G(\tilde{K})_{E_1, E_2, E_3})^{\circ}$ 
for $a \in \mathcal{S}_1(1, \tilde{K})^{\circ}$. 
And 
$\delta_i(-a) = \sigma_i \delta_i(a) 
= \delta_i(a) \sigma_i$ 
with 
$\sigma_i \in 
((G(\tilde{K})_{E_i})^{\circ})_{E_{i+1}, E_{i+2}}$ 
by (1)(\roman{f1}). 
By (\roman{f1}), 
$\{ \delta_i(a) |~a \in \mathcal{S}_1(1, \tilde{K}) \} 
= \{ \delta_i(a), \delta_i(-a) 
|~a \in \mathcal{S}_1(1, \tilde{K})^{\circ} \} 
= \{ \delta_i(a), \delta_i(a) \sigma_i 
|~a \in \mathcal{S}_1(1, \tilde{K})^{\circ} \} 
\subseteqq  
(G(\tilde{K})_{E_i})^{\circ})_{E_{i+2}, E_{i+3}}$. 
By the last claim of (\roman{f1}), 
$\{ \delta_i(a) |~a \in \mathcal{S}_1(1, \tilde{K}) \} 
= \{ \delta_i(a) |~a \in \mathcal{S}_1(1, \tilde{K})^{\circ} \} 
\subseteqq (G(\tilde{K})_{E_1, E_2, E_3})^{\circ}$ 
when 
$\tilde{K} = \bH', \bO'; \bC^{\mathbb{C}}, 
\bH^{\mathbb{C}}, \bO^{\mathbb{C}}$.

(\roman{f3}) 
$\tilde{K}$ is associative 
by $d_{\tilde{K}} \leqq 4$. 
Hence, 
$GL_{\mathbb{F}}(\mathcal{J}_3(\tilde{K}))
_{E_1, E_2, E_3, F_i(1)} \ni \beta_i(a)$ 
is well-defined such that 
$\beta_i(a)^{-1}= \delta_i(\bar{a})$, 
$\beta_i(a) \tau = \tau \beta_i(\tau a)$, 
$\beta_i(a) E = E$. 
Then 
$\mathrm{det}(\beta_i(a) \mathbb{X}(r;x)) 
= r_1 r_2 r_3 
+ 2 ( \overline{a x_i \overline{a}} | 
(a x_{i+1}) (x_{i+2} \overline{a}) ) 
- r_i N(a x_i \overline{a}) 
- r_{i+1} N(a x_{i+1}) 
- r_{i+2} N(x_{i+2} \overline{a}) 
= r_1 r_2 r_3 
+ (2 (\bar{x}_i | x_{i+1} x_{i+2}) 
- r_i N(x_i)) N(a)^2 
- (r_{i+1} N(x_{i+1}) 
- r_{i+2} N(x_{i+2})) N(a) 
= \mathrm{det}(\mathbb{X}(r; x))$ 
by Lemma 1.1 (1). 
Because of Proposition 0.1 (1), 
$\beta_i(a) \in 
(G(\tilde{K})_{E_1, E_2, E_3, F_1(1)})^{\circ}$ 
for $a \in \mathcal{S}_1(1, \tilde{K})^{\circ}$. 
By (1) (\roman{f1}), 
$\sigma_i \in 
((G(\tilde{K})_{E_i})^{\circ})_{E_{i+1}, E_{i+2}}$ 
with 
$\beta_i(-a) = \sigma_i \beta_i(a) 
= \beta_i(a) \sigma_i$.  
By virtue of (\roman{f1}), 
$\{ \beta_i(a) |~a \in \mathcal{S}_1(1, \tilde{K}) \} 
= \{ \beta_i(a), \beta_i(-a) 
|~a \in \mathcal{S}_1(1, \tilde{K})^{\circ} \} 
= \{ \beta_i(a), \beta_i(a) \sigma_i 
|~a \in \mathcal{S}_1(1, \tilde{K})^{\circ} \}$ 
is contained in 
$(G(\tilde{K})_{E_i})^{\circ})_{E_{i+2}, E_{i+3}, F_i(1)}$.  
By the last claim of (\roman{f1}), 
$\{ \beta_i(a) |~a \in \mathcal{S}_1(1, \tilde{K}) \} 
= \{ \beta_i(a) |~a \in 
\mathcal{S}_1(1, \tilde{K})^{\circ} \}$ 
is contained in 
$(G(\tilde{K})_{E_1, E_2, E_3, F_i(1)})^{\circ}$ 
when 
$\tilde{K} = \bH'; \bC^{\mathbb{C}}, 
\bH^{\mathbb{C}}$ with 
$d_{\tilde{K}} \leqq 4$. 
\qed

\bigskip

Put 
$G_J(\tilde{K}_{\tau}) 
:= \{ \beta_j(t; a, \sqrt{-1}) 
|~j \in J, t \in \mathbb{R}, 
a \in \tilde{K}_{\tau}, 
N(a) = 1 \}$ 
for any subset 
$J \subseteqq  \{ 1, 2, 3 \}$. 
By Lemma 1.4 (1) (\roman{f1}), 
$G_J(\tilde{K}_{\tau}) 
\subset (G(\tilde{K})^{\tau})^{\circ}$. 
By Proposition 0.1 (2), 
$G(\tilde{K})^{\tau}$ 
and the identity connected component 
$(G(\tilde{K})^{\tau})^{\circ}$ 
are compact.

\bigskip

{\sc Lemma 1.5.} 
(1) 
{\it 
For any $X \in \mathcal{J}_3(\tilde{K})_{\tau}$ 
and any closed subgroup 
$H$ of $G(\tilde{K})^{\tau}$ 
such that 
$G_J(\tilde{K}_{\tau}) \subseteqq H$ 
with some $J \subseteqq \{ 1, 2, 3 \}$, 

\[
{\cal O}_H(X) \cap 
\{ Y \in \mathcal{J}_3(\tilde{K}) 
|~(Y | F_j(x)) 
= 0~(j \in J,~x \in \tilde{K}_{\tau}) \} 
\ne \emptyset. 
\]
}

\smallskip

(2) 
{\it 
${\cal O}_{(G(\tilde{K})^{\tau})^{\circ}}(X) 
\cap \{ \mathrm{diag}(r_1, r_2, r_3) 
|~ r_i \in \mathbb{R}~(i = 1, 2, 3) \} 
\ne \emptyset$ 
for any $X \in \mathcal{J}_3(\tilde{K})_{\tau}$, 
where 
$\{ r_1, r_2, r_3 \} = \Lambda_X$ 
iff 
$\mathrm{diag}(r_1, r_2, r_3) \in 
{\cal O}_{(G(\tilde{K})^{\tau})^{\circ}}(X)$. 
}

\smallskip

(3) 
{\it 
${\cal O}_{G(K)^{\circ}}(X) 
\cap \{ Y + \sqrt{-1} \mathrm{diag}(r_1, r_2, r_3) 
|~Y \in \mathcal{J}_3(K),
~r_i \in \mathbb{R}~(i = 1, 2, 3) \} 
\ne \emptyset$ 
for any $X \in \mathcal{J}_3(K^{\mathbb{C}})$.
}

\smallskip

(4) 
{\it 
${\cal O}_{(G(K')^{\tau})^{\circ}}(X) 
\cap \{ \mathbb{X}(s; y) 
|~s_i \in \mathbb{R},
~y_i = \sqrt{-1} p_i e_4,
~p_i \in K \cap K'~(i = 1, 2, 3) \} 
\ne \emptyset$ 
for any $X \in \mathcal{J}_3(K')$. 
}

\medskip

{\it Proof}. 
(1) (cf. \cite[3.3]{Yi1968}): 
Since the closed subgroup $H$ of 
the compact group 
$G(\tilde{K})^{\tau}$ 
is compact, 
the orbit 
$\mathcal{O}_{H}(X)$ is compact, 
which is contained in 
$\mathcal{J}_3(\tilde{K})_{\tau}$ 
if $X \in \mathcal{J}_3(\tilde{K})_{\tau}$. 
Put 
$\phi: 
\mathcal{J}_3(\tilde{K})_{\tau} \longrightarrow \mathbb{R}; 
\mathbb{X}(r; x) \mapsto \sum_{j=1}^3 r_j^2$, 
which is a continuous $\mathbb{R}$-valued function 
admitting a maximal point 
$\mathbb{X}(r; x) \in \mathcal{O}_H(X)$. 
Suppose that 
$(\mathbb{X}(r, x) | F_j(q)) \ne 0$ 
for some $j \in J$ and $q \in \tilde{K}_{\tau}$. 
By Lemma 1.1 (1), $2 (x_j| q) \ne 0$. 
Since $(x| y)$ is non-degenerate on 
$\tilde{K}$ and $\tilde{K}_{\tau}$, 
$\tilde{K} 
= \tilde{K}_{\tau} \oplus \tilde{K}_{\tau}^{\perp}$ 
for 
$\tilde{K}_{\tau}^{\perp} 
:= \{ x \in \tilde{K} 
|~(x| y) = 0~(y \in \tilde{K}_{\tau}) \}$, 
so that $x_j = y_j + y_j^{\perp}$ for some 
$y_j \in \tilde{K}_{\tau}$ 
and $y_j^{\perp} \in \tilde{K}_{\tau}^{\perp}$. 
Then $(y_j| q) = (x_j| q) \ne 0$, 
so that 
$y_j \ne 0$ 
and  
$(y_j| y_j) = (\tau y_j| y_j) > 0$. 
Put 
$a := y_j / \sqrt{(y_j|y_j)} \in \tilde{K}_{\tau}$, 
so that $(a| a) = 1$. 
By Lemma 1.4 (1) (\roman{f1}), 
$\beta_j(t; a, \sqrt{-1}) 
\in G_J(\tilde{K}_{\tau}) \subseteqq H$. 
Put 
$\varepsilon := (a| x_j) = (a| y_j) 
= \sqrt{(y_j|y_j)} > 0$, 
$s_j^{\pm} := (r_{j+1} \pm r_{j+2})/2 \in \mathbb{R}$ 
and 
$Y(t) := \beta_j(t; a, \sqrt{-1}) \mathbb{X}(r; x) 
\in \mathcal{J}_3(\tilde{K})_{\tau}$. 
Then 
$\phi(Y(t)) 
= r_j^2 + \sum_{\pm} 
( s_j^{+} 
\pm ( s_j^{-} \cos(2t) + \varepsilon \sin (2t) ) )^2 
= r_j^2 + 2 ( s_j^{+} )^2 
+ 2 ( s_j^{-} \cos(2t) + \varepsilon \sin (2t) )^2 
= r_j^2 + 2 ( s_j^{+} )^2 
+ 2 ( (s_j^{-} )^2 + \varepsilon^2 ) \cos(2 t + \theta)$ 
for some constant $\theta$ of $t$ 
determined by $s_j^{-}$ and $\varepsilon > 0$. 
Hence, 
$\phi(Y(\frac{-\theta}{2})) 
= r_j^2 + r_{j+1}^2 + r_{j+2}^2 + 2 \varepsilon^2 
= \phi(\mathbb{X}(r; x)) + 2 \varepsilon^2 
\geqq \phi(Y(\frac{-\theta}{2})) + 2 \varepsilon^2$ 
by the maximality of $\phi(\mathbb{X}(r; x))$, 
which gives $\varepsilon = 0$, a contradiction.

(2) 
Take 
$X \in \mathcal{J}_3(\tilde{K})_{\tau}$. 
By (1) on 
$H := (G(\tilde{K})^{\tau})^{\circ} 
\supseteqq G_{\{ 1, 2, 3 \}}(\tilde{K}_{\tau})$, 
there exists $\beta \in H$ 
such that $(\beta X, F_i(x)) = 0$ 
($x \in \tilde{K}_{\tau}; i = 1, 2, 3$), 
so that 
$\beta X = \mathrm{diag}(r_1, r_2, r_3)$ 
for some $r_i \in \mathbb{R}$ 
($i = 1, 2, 3$). 
In this case, 
$\Phi_X(\lambda) 
= \Phi_{\mathrm{diag}(r_1, r_2, r_3)}(\lambda) 
= \Pi_{i = 1}^3 (\lambda - r_i)$, 
so that 
$\{ r_1, r_2, r_3 \} = \Lambda_X$. 
Conversely if 
$\{ r_1, r_2, r_3 \} = \Lambda_X$, 
then 
$\mathrm{diag}(s_1, s_2, s_3) \in 
{\cal O}_{(G(\tilde{K})^{\tau})^{\circ}}(X)$ 
for some $\{ s_1, s_2, s_3 \} = \Lambda_X$, 
so that 
$\mathrm{diag}(r_1, r_2, r_3) \in 
{\cal O}_{(G(\tilde{K})^{\tau})^{\circ}}(X)$ 
by Lemma 1.4 (1) (\roman{f2}).

(3) 
Take $X \in \mathcal{J}_3(K^{\mathbb{C}})$. 
Then $X = X_1 + \sqrt{-1} X_2$ 
for some 
$X_i \in \mathcal{J}_3(K) 
= \mathcal{J}_3(K^{\mathbb{C}})_{\tau}$ 
($i = 1, 2$). 
By (2), 
there exist 
$\beta \in 
(G(K^{\mathbb{C}})^{\tau})^{\circ} = G(K)^{\circ}$ 
and $\{ r_1, r_2, r_3 \} \subset \mathbb{R}$ 
such that 
$\beta X_2 = \mathrm{diag}(r_1, r_2, r_3)$, 
so that 
$\beta X = \beta X_1 + \sqrt{-1} \beta X_2$ 
has the required form with  
$\beta X_1 \in \mathcal{J}_3(K)$.

(4) 
Take $X \in \mathcal{J}_3(K')$. 
Then $X = X_+ + X_-$ 
for some 
$X_{\pm} \in 
\mathcal{J}_3(K')_{\pm \tau}$. 
By (1) on 
$H := (G(K')^{\tau})^{\circ} 
\supseteqq 
G_{\{1, 2, 3 \}}(\tilde{K}_{\tau})$, 
there exists $\beta \in H$ 
such that 
$\beta X_+ = \mathrm{diag}(r_1, r_2, r_3)$ 
for some $r_i \in \mathbb{R}$. 
Then 
$\beta X = \beta X_+ + \beta X_-$ 
has the required form because of 
$\beta X_- \in \mathcal{J}_3(K')_{-\tau} 
= \{ \mathbb{X}(0; y) 
|~y_i = \sqrt{-1} p_i e_4;~p_i \in K \cap K'
~(i = 1, 2, 3) \}$.~\qed

\bigskip

{\sc Lemma 1.6.} 
(1) 
{\it 
For a positive integer $m$, 
let $f(X_1, \cdots, X_m)$ be 
a $\mathcal{J}_3(\tilde{K})$-valued 
polynomial of $E$ and 
$X_1, \cdots X_m \in \mathcal{J}_3(\tilde{K})$ 
with respect to $\circ$, $\times$ 
and the scalar multiples of 
$\mathrm{tr}(X_i)$, $(X_i | X_j)$, 
$\mathrm{det}(X_i)$ and $(X_i | X_j | X_k)$ 
for $i, j, k \in \{ 1, \cdots, m \}$. 
Assume that 
$f(X_1, X_2, \cdots, X_m) = 0$ 
for any 
$X_2, \cdots X_m \in \mathcal{J}_3(K)$ 
and all diagonal forms 
$X_1$ in $\mathcal{J}_3(\mathbb{R})$.  
Then \\
$f(X_1, \cdots, X_m) = 0$ 
for all 
$X_1, \cdots X_m \in \mathcal{J}_3(K^{{\mathbb{C}}})$. 
}

\smallskip
\noindent
(2) 
{\it 
Assume that 
$X, Y \in \mathcal{J}_3(\tilde{K})$. 
Then:
} 

\smallskip

(\roman{f1}) $X \circ ( (X \circ X) \circ Y ) 
= (X \circ X) \circ (X \circ Y)$; 

\smallskip

(\roman{f2}) 
$X^{\times 2} \circ X = \mathrm{det}(X)E$,
~$(X^{\times 2})^{\times 2} = \mathrm{det}(X) X$; 

\smallskip

(\roman{f3}) 
$X^{\times 2} \times X 
= - \frac{1}{2} 
\{ 
\mathrm{tr}(X) X^{\times 2} 
+ \mathrm{tr}(X^{\times 2}) X 
- (\mathrm{tr}(X^{\times 2}) \mathrm{tr}(X) 
- \mathrm{det}(X))E 
\}$. 

\smallskip
\noindent
(3) 
{\it 
$V_X$ is the minimal subspace over 
$\mathbb{F}$ generated by $X$ and $E$ 
under the cross product. 
Especially, 
$\varphi_X (\lambda)^{\times 2} \in V_X$ 
for all $\lambda \in \mathbb{F}$. 
}

\smallskip
\noindent
(4) 
$(\varphi_X (\lambda_1)^{\times 2})^{\times 2} 
= 0$ 
{\it 
if $X \in \mathcal{J}_3(\tilde{K})$ 
and $\lambda_1 \in {\mathbb{C}}$ 
with $\Phi_X(\lambda_1) = 0$. 
}

\smallskip
\noindent 
(5) 
{\it 
$\mathcal{M}_{2 3}(\tilde{K}) 
= \{ X \in \mathcal{J}_3(\tilde{K})_0 
|~X^{ \times 2} \ne 0,~\mathrm{tr}(X^{\times 2}) = 0,
~(X^{\times 2})^{\times 2} = 0 \}$ 
and 
$\{ X^{\times 2} |~X \in \mathcal{M}_{2 3}(\tilde{K}) \} 
\subseteqq \mathcal{M}_1(\tilde{K})$. 
}

\medskip

{\it Proof}. 
(1) (cf. \cite[p.42]{Fh1951}, \cite[p.74, $\ell \ell$.2--4]{Yi1975}, 
\cite[p.91, Corollary \Roman{f5}.2.6]{FK1994}): 
By Lemma 1.5 (2), 
any $X_1 \in \mathcal{J}_3(K)$ 
admits some 
$\beta \in G(K) = G(K^{\mathbb{C}})^{\tau}$ 
such that $\beta X_1$ 
is a diagonal form in $\mathcal{J}_3(\mathbb{R})$.  
Then 
$f(\beta X_1, X_2, \cdots, X_m) = 0$ 
for any  
$X_i \in \mathcal{J}_3(K)$ 
with $i \in \{ 2, \cdots, m \}$. 
By Proposition 0.1, 
$\beta$ preserves 
$\circ$, $\times$, $\mathrm{tr}(*)$, 
$(* | *)$, $\mathrm{det}(*)$, $(* | * | *)$ and $E$, 
so that  
$f(X_1, \beta^{-1} X_2 \cdots, \beta^{-1} X_m) = 0$ 
for all 
$X_i \in \mathcal{J}_3(K)~(i = 2, \cdots, m)$. 
Hence, 
$f(X_1, \cdots, X_m) = 0$ 
for all 
$X_i \in \mathcal{J}_3(K)~(i = 1, \cdots, m)$. 
Since this formula consists of some polynomial equations 
on the $\mathbb{R}$-coefficients of 
each matrix entry of $X_i$'s 
with respect to the $\mathbb{R}$-basis $\{ e_j \}$ of $K$, 
the formula holds on $\mathcal{J}_3(K^{{\mathbb{C}}}) = 
\mathcal{J}_3(K) \otimes_{\mathbb{R}} {\mathbb{C}}$. 

(2) 
The formulas in (\roman{f1}) and (\roman{f2}) 
are polynomials of $X, Y$ and $E$ 
with respect to 
$\circ$, $\times$, $\mathrm{tr}(*)$, 
$(* | *)$, $\mathrm{det}(*)$, $(* | * | *)$. 
If $X$ is a diagonal form in ${\cal J}_3(\mathbb{R})$, 
the formulas can be checked by Lemma 1.1 (1), 
easily. 
By (1), the formulas (\roman{f1}) and (\roman{f2}) 
hold for any $X, Y \in {\cal J}_3(K^{\mathbb{C}})$. 
Hence, they hold for any $X, Y \in {\cal J}_3(\tilde{K}) 
\subseteqq {\cal J}_3(K^{\mathbb{C}})$. 
The formula (\roman{f3}) follows from 
the first formula of (\roman{f2}) 
and the definition of cross product with 
$(X^{\times 2} | X) = 3~\mathrm{det}(X)$.

(3) follows from the formulas 
in (\roman{f2}), (\roman{f3}) and Lemma 1.1 (3).

(4) 
$(\varphi_X (\lambda_1)^{\times 2})^{\times 2}  
= \mathrm{det}(\varphi_X (\lambda_1)) \varphi_X (\lambda_1) 
= \Phi_X(\lambda_1) \varphi_X(\lambda_1) = 0$ 
by the second formula of (\roman{f2}) in (2).

(5) 
By (2) (\roman{f2}), 
$(X^{\times 2})^{\times 2} = \mathrm{det}(X) X$, 
so that $(X^{\times 2})^{\times 2} = 0$ 
if and only if $\mathrm{det}(X) = 0$, 
which gives the results.~\qed
 
\bigskip

The formula (\roman{f1}) 
of Lemma 1.6 (2) implies that  
$(\mathcal{J}_3(\tilde{K}), \circ)$ 
is a Jordan algebra over $\mathbb{F}$, 
which is also reduced simple 
in the sense of N.~Jacobson 
\cite[Chapters \Roman{f4}, \Roman{f9}]{Jn1968}, 
where $\mathcal{J}_3(\tilde{K})$ 
is called {\it split} iff 
$\tilde{K}$ is {\it split} (i.e.  non-division), 
that is the case when $\tilde{K} = K'$ 
or ${K'}^{{\mathbb{C}}}$. 

\bigskip

{\it Proof of ``Proposition 0.1 (3) when $\tilde{K} = K$''} 
(cf, \cite{Fh1951}, \cite[4.1 Proposition]{Yi1968}). 
Take any 
$X \in \mathcal{P}_2(K) 
\subset \mathcal{J}_3(K) 
= \mathcal{J}_3(K^{\mathbb{C}})^{\tau}$. 
By Lemma 1.5 (2) with $\tilde{K} = K^{\mathbb{C}}$, 
there exists 
$\alpha \in G(K)^{\circ} 
= (G(K^{\mathbb{C}})^{\tau})^{\circ}$ 
such that 
$\alpha X = \mathrm{diag}(r_1, r_2, r_3)$ 
for some $r_1, r_2, r_3 \in \bR$. 
By $\mathrm{tr}(\alpha X) = 1$ 
and $(\alpha X)^{\times 2} = 0$, 
$r_1 + r_2 + r_3 = 1$ 
and $r_2 r_3 = r_3 r_1 = r_1 r_2 = 0$, 
so that 
$(r_1, r_2, r_3) = (1, 0, 0), (0, 1, 0), (0, 0, 1)$. 
By Lemma 1.4 (1) (\roman{f2}), there exists 
$\hat{\beta} \in (G(K^{\mathbb{C}})^{\tau})^{\circ} 
= G(K)^{\circ}$ such that 
$\hat{\beta} (\alpha X) = E_1$.
\qed

\bigskip

H.~Freudenthal \cite[5.1]{Fh1951} 
gave the diagonalization theorem on  
$\mathcal{J}_3(\bO)$ with the action of 
$\{ \alpha \in G(\bO) 
|~\mathrm{tr}(\alpha X) = \mathrm{tr}(X) \}$ 
(cf. \cite[3.3 Theorem]{Yi1968}, 
\cite[p.206, Lemma 1]{Ms1989}, 
\cite[Proposition 1.4]{So2006}, 
\cite[p.90, Theorem \Roman{f5}.2.5]{FK1994}), 
which is developing to Lemma 1.5 (2) for 
$\tilde{K} = \bO^{\mathbb{C}}$ 
with ${\cal J}_3(\tilde{K})_{\tau} 
= {\cal J}_3(\bO^{\mathbb{C}})_{\tau} 
= {\cal J}_3(\bO)$ 
under the action of 
$G(\tilde{K})^{\tau} 
= G(\bO^{\mathbb{C}})^{\tau} \cong G(\bO) =: F_4$. 
I.~Yokota \cite[4.2 and 6.4 Theorems]{Yi1968} 
proved the connectedness and the simply connectedness of 
$F_4$ by the diagonalization theorem of H. Freudenthal 
(cf. \cite[Appendix]{My1954}, 
\cite[p.210, Theorem 3]{Ms1989}, 
\cite[p.175, Proposition 1.4]{GOV1994}). 
O.~Shukuzawa \& 
I.~Yokota \cite[p.3, Remark]{SY1979} 
(cf. \cite[p.63, Theorem 9; p.54, Remark]{Yi1977}) 
proved the connectedness of $F_4' := G(\bO')$ 
by showing the first formula of Proposition 0.1 (1) 
by virtue of Hamilton-Cayley formula on 
$\mathcal{J}_3(\bO')$ 
given as the first formula of Lemma 1.6(2)(\roman{f2}) 
(cf. \cite[p.119, Proposition 5.1.5]{SV}, 
\cite[Lemma 14.96]{Hfr1990}). 
Because of 
$F_4^{\mathbb{C}} 
\cong (F_4^{\mathbb{C}})^{\tau} \times \mathbb{R}^{52}$ 
with 
$(F_4^{{\mathbb{C}}})^{\tau} = F_4$ 
\cite[Theorem 2.2.2]{Yi1990} (cf. Proposition 0.1 (2)), 
$F_4^{\mathbb{C}} := G(\bO^{\mathbb{C}})$ 
is connected and simply connected, 
so that 
$F_4' = (F_4^{{\mathbb{C}}})^{\tau \gamma}$ 
is again proved to be connected 
by virtue of a theorem of P.K.~Rasevskii \cite{Rpk1974}. 

\bigskip

\begin{center}
{\bf 2. 
Proposition 0.1 (3) and (4) (\roman{f1}).} 
\end{center} 

Assume that 
$\tilde{K} = K'$ or $K^{\mathbb{C}}$ 
with 
$K' = \bC', \bH'$ or $\bO'$; 
and 
$K^{\mathbb{C}} = \mathbb{R}^{\mathbb{C}}, 
\bC^{\mathbb{C}}, \bH^{\mathbb{C}}$ 
or $\bO^{\mathbb{C}}$. 
And put $\sigma := \sigma_1$ 
defined in Lemma 1.4 (1) (\roman{f1}) such that 
$\sigma^2 = \mathrm{id}_{\mathcal{J}_3(\tilde{K})}$. 
Then 
$\mathcal{J}_3(\tilde{K}) = 
\mathcal{J}_3(\tilde{K})_{\sigma} 
\oplus 
\mathcal{J}_3(\tilde{K})_{-\sigma}$, 
$\mathcal{J}_3(\tilde{K})_{\sigma} 
= \{ \sum_{i=1}^3 r_i E_i + F_1(x_1) 
|~r_i \in \mathbb{F}, x_1 \in \tilde{K} \}$ 
such that 
$\mathcal{J}_3(\tilde{K})_{-\sigma} 
= \{ F_2(x_2) + F_3(x_3) |~x_2, x_3 \in \tilde{K} \} 
= \{ X \in \mathcal{J}_3(\tilde{K}) 
|~(X, Y) = 0
~(Y \in \mathcal{J}_3(\tilde{K})_{\sigma}) 
\}$. 
And 
$\mathcal{J}_3(\tilde{K})_{\sigma} 
= \mathbb{F} E_1 \oplus \mathcal{J}_2(\tilde{K})$ 
with 
$\mathcal{J}_2(\tilde{K}) 
:= \{ \sum_{i = 2}^3 r_i E_i + F_1(x_1) 
|~r_i \in \mathbb{F}, 
x_1 \in \tilde{K} \}$. 
By Lemma 1.1 (1), 
$\mathcal{J}_3(\tilde{K})_{L^{\times}_{2E_1}} 
=\{ r ( E_2 + E_3 ) |~r \in \mathbb{F} \}$ 
and 
$\mathcal{J}_3(\tilde{K})_{- L^{\times}_{2E_1}} 
= \{ r (E_2 - E_3) + F_1(x) 
|~r \in \mathbb{F}, x \in \tilde{K} \}$, 
so that 
$\mathcal{J}_2(\tilde{K}) 
= \mathcal{J}_3(\tilde{K})_{L^{\times}_{2E_1}} 
\oplus \mathcal{J}_3(\tilde{K})_{- L^{\times}_{2E_1}}$.

\bigskip

{\sc Lemma 2.1}. 
(1) 
{\it 
$G(\tilde{K})_{E_1} 
= G(\tilde{K})_{E_1, E_2 + E_3, 
\mathcal{J}_3(\tilde{K})_{\pm L^{\times}_{2 E_1}}, 
\mathcal{J}_2(\tilde{K}), 
\mathcal{J}_3(\tilde{K})_{\pm \sigma}}$. 
} 

\smallskip

(2) 
(\roman{f1}) 
{\it 
$\{ \mathbb{X}((X | E_1), s_2, s_3; 0, 0, 0) 
|~s_2, s_3 \in \mathbb{R}; s_2 \geqq s_3 \} 
\cap {\cal O}_{(G(K)_{E_1})^{\circ}}(X) 
\ne \emptyset$ 
and 
$\{ \mathbb{X}((X | E_1), t_2, t_3; u, 0, 0) 
|~t_2, t_3, u \in \mathbb{R}; u \geqq 0 \} 
\cap {\cal O}_{((G(K)_{E_3})^{\circ})_{E_1, E_2}}(X) 
\ne \emptyset$ 
if $X \in \mathcal{J}_3(K)_{\sigma}$.
}

\smallskip

(\roman{f2}) 
{\it 
$\{ \mathbb{X}((X | E_1), s_2, s_3; u \sqrt{-1} e_4, 0, 0) 
|~s_2, s_3, u \in \mathbb{R}; u \geqq 0, s_2 \geqq s_3 \} 
~\cap {\cal O}_{(G(K')^{\circ})^{\tau}_{E_1}}(X) 
\ne \emptyset$ 
if $X \in \mathcal{J}_3(K')_{\sigma}$.}

\smallskip

(\roman{f3}) 
{\it 
$\{ \mathbb{X}((X | E_1), t_2 + \sqrt{-1} s_2, 
t_3 + \sqrt{-1} s_3; u, 0, 0) 
|~t_2, t_3, s_2, s_3, u \in \mathbb{R}; 
u \geqq 0, s_2 \geqq s_3 \} 
\cap {\cal O}_{(G(K)^{\circ})_{E_1}}(X) 
\ne \emptyset$ 
if $X \in \mathcal{J}_3(K^{\mathbb{C}})_{\sigma}$.}

\medskip

{\it Proof}. 
(1) 
For $\alpha \in G(\tilde{K})_{E_1}$, 
one has that 
$\alpha (E_2 + E_3) = \alpha (E - E_1) 
= \alpha E - \alpha E_1 
= E - E_1 = E_2 + E_3$ and 
$\alpha \mathcal{J}_3(\tilde{K})_{\pm L^{\times}_{2 E_1}} 
= \mathcal{J}_3(\tilde{K})_{\pm L^{\times}_{2 E_1}}$ 
since $\alpha$ preserves $\times$ by Proposition 0.1 (1), 
so that 
$\alpha \mathcal{J}_2(\tilde{K}) 
= \mathcal{J}_2(\tilde{K})$, 
$\alpha \mathcal{J}_3(\tilde{K})_{\sigma} 
= \mathcal{J}_3(\tilde{K})_{\sigma}$ 
and 
$\alpha \mathcal{J}_3(\tilde{K})_{- \sigma} 
= \mathcal{J}_3(\tilde{K})_{- \sigma}$ 
because of the orthogonal direct-sum decompositions of them 
and that $\alpha$ preserves $(* | *)$ on 
$\mathcal{J}_3(\tilde{K})$ by Proposition 0.1 (1).

(2) (\roman{f1}) 
Take any $X \in \mathcal{J}_3(K)_{\sigma}$. 
Then there exist 
$r_i \in \mathbb{R}$ and $x_1 \in K$ 
such that 
$X = \mathbb{X}(r_1, r_2, r_3; x_1, 0, 0)$. 
By Lemmas 1.4 (1) and 1.5 (1) 
with $\tilde{K} = K = K_{\tau}$, 
$F = \mathbb{R}$ and 
$H := (G(K)_{E_1})^{\circ} 
\supseteqq G_J (K)$ 
with $J = \{ 1 \}$, 
there exists 
$\alpha \in H$ such that 
$(\alpha X)_{F_1} = 0$. 
By (1), 
$\alpha X \in \mathcal{J}_3(K)_{\sigma}$, 
so that $\alpha X$ is diagonal with 
$s_i := (\alpha X | E_i) \in \mathbb{R}$ 
($i = 1, 2, 3$) such that 
$s_1 = (\alpha X | \alpha E_1) 
= (X | E_1) = r_1$. 
If $s_2 \geqq s_3$, then 
$\alpha X$ gives an element of 
the left-handed set of the first formula. 
If $s_2 < s_3$, 
put $\alpha_1 := \hat{\beta}_1 \alpha$ 
with $\hat{\beta}_1 \in (G(K)_{E_1})^{\circ}$ 
given in Lemma 1.4 (1) (\roman{f2}), 
so that $\alpha_1 X$ gives an element of 
the left-handed set of the first formula. 
Hence, follows the first formula.

If $x_1 = 0$, then $X$ gives an element of 
the left-handed side of the second formula 
with $u = 0 \in \mathbb{R}$. 
If $x_1 \ne 0$, put 
$a := x_1 / \sqrt{(x_1| x_1)} \in \mathcal{S}_1(1, K)$, 
so that 
$\delta_3(a) \in ((G(K)_{E_3})^{\circ})_{E_1, E_2}$ 
in Lemma 1.4 (2) (\roman{f2}) 
such that 
$\delta_3(a) X = 
\mathbb{X}(r_1, r_2, r_3; u, 0, 0)$ 
with 
$u := \sqrt{(x_1 | x_1)}> 0$, 
which gives an element of the left-handed side 
of the second formula. 
Hence, follows the second formula.

(\roman{f2}) 
Take any 
$X \in \mathcal{J}_3(K')_{\sigma}$. 
By Lemmas 1.4 (1) and 1.5 (1) 
with $\tilde{K} = K'$ and 
$H := (G(K')^{\tau}_{E_1})^{\circ} 
\supseteqq G_J (K_{\tau}')$ 
with $J = \{1 \}$, there exists 
$\beta \in H$ such that 
$(\beta X | F_1(x)) = 0$ 
for all $x \in K_{\tau}' = K \cap K'$. 
By (1), 
$\beta X \in \mathcal{J}_3(K')_{\sigma}$. 
Hence, 
$\beta X = 
\mathbb{X}((X|E_1), s_2, s_3; \sqrt{-1} q e_4, 0, 0)$ 
for some $q \in K \cap K'$. 
Put $\alpha_1 := \beta$ (if $s_2 \geqq s_3$) 
or $\hat{\beta}_1 \beta$ (if $s_2 < s_3$), 
so that 
$\alpha_1 \in H$ by Lemma 1.4 (1) (\roman{f2}). 
Then 
$\alpha_1 X  
= \mathbb{X}((X|E_1), s_2, s_3; \sqrt{-1} q e_4, 0, 0)$ 
for some $q \in K \cap K'$, $s_2, s_3 \in \mathbb{R}$ 
with $s_2 \geqq s_3$. 
Put $\alpha := \alpha_1$ (if $q = 0$) 
or $\delta_3(a) \beta$ 
for $a := q /\sqrt{(q| q)} \in K_{\tau}'$ 
with $N(a) = 1$ (if $q \ne 0$), 
where 
$\delta_3(a) \in 
((G(K')_{E_3})^{\circ})^{\tau}_{E_1, E_2} 
\subseteqq (G(K')^{\circ})^{\tau}_{E_1}$ 
by Lemma 1.4 (2) (\roman{f2}). 
Then 
$\alpha X = \mathbb{X}((X | E_1), s_2, s_3; 
\sqrt{-1} u e_4, 0, 0)$ 
with $u := \sqrt{N(q)} \geqq 0$, 
which is an element of the left-handed set.

(\roman{f3}) 
Take any 
$X \in \mathcal{J}_3(K^{\mathbb{C}})_{\sigma}$. 
Then 
$X = X_1 + \sqrt{-1} X_2$ for some 
$X_i \in \mathcal{J}_3(K)_{\sigma}$ 
($i = 1, 2$). 
By (\roman{f1}), 
there exist 
$\alpha_1 \in (G(K)_{E_1})^{\circ}$ 
such that 
$\alpha_1 (X_2) = \mathbb{X}((X_2 | E_1), s_2, s_3; 0, 0, 0)$ 
for some $s_2, s_3 \in \mathbb{R}$ with $s_2 \geqq s_3$. 
Because of 
$\mathcal{J}_3(K)_{\sigma} 
\ni \alpha_1 (X_1) = 
\mathbb{X}((X_1 | E_1), t_2, t_3; x, 0, 0)$ 
for some 
$t_2, t_3 \in \mathbb{R}$ and $x \in K$, 
so that 
$\alpha_1(X) 
= \mathbb{X}((X | E_1), 
t_2 + \sqrt{-1} s_2, 
t_3 + \sqrt{-1} s_3; 
x, 0, 0)$. 
Put 
$\alpha := \alpha_1$ (if $x = 0$) 
or $\delta_3(a) \alpha_1$ 
with $a := x/\sqrt{(x| x)} \in \mathcal{S}_1(1, K)$
(if $x \ne 0$), 
where 
$\delta_3(a) \in ((G(K)_{E_3})^{\circ})_{E_1, E_2}$ 
by Lemma 1.4 (2) (\roman{f2}). 
Then 
$\alpha \in (G(K)^{\circ})_{E_1}$ 
and 
$\alpha X 
= \mathbb{X}((X | E_1), 
t_2 + \sqrt{-1} s_2, 
t_3 + \sqrt{-1} s_3; 
u, 0, 0)$ with 
$u := \sqrt{(x| x)} \geqq 0$, 
which is an element of the left-handed set.~\qed

\bigskip

For $c \in \mathbb{F}$, put 
$\mathcal{S}_2(c, \tilde{K}) 
:= \{ W \in \mathcal{J}_3(\tilde{K})_{- L^{\times}_{2E_1}} 
|~(W | W) = c, W \ne 0 \}$, 
which is said to be 
{\it a generalized sphere 
of second kind over} 
$\mathbb{F}$. 
Then 
$G(\tilde{K})_{E_1} 
= \cap_{c \in \mathbb{F}}
~G(\tilde{K})_{E_1, \mathcal{S}_2(c, \tilde{K}).}$  

\bigskip

{\sc Lemma 2.2.} 
(1) 
{\it 
$\mathcal{J}_3(K')_{-L^{\times}_{2E_1}} 
= (\cup_{c \in \mathbb{R}} \mathcal{S}_2(c, K')) 
\cup \{ 0 \}$ 
such that   
}

\smallskip

\noindent 
(\roman{f1}-1) 
{\it 
$\mathcal{S}_2(c, K') 
= \mathcal{O}_{(G(K')^{\circ})_{E_1}} 
(\sqrt{\frac{c}{2}} (E_2 - E_3))$ 
for $c > 0$;  
}

\smallskip

\noindent 
(\roman{f1}-2) 
{\it 
$\mathcal{S}_2(c, K') 
= \mathcal{O}_{(G(K')^{\circ})_{E_1}}
(\sqrt{\frac{- c}{2}} F_1(\sqrt{-1} e_4))$ 
for $c < 0$; 
}

\smallskip

\noindent 
(\roman{f2}) 
{\it 
$\mathcal{S}_2(0, K') 
= \mathcal{O}_{(G(K')^{\circ})_{E_1}}(M_{1'})$; 
and 
}

\smallskip

\noindent 
(\roman{f3}) 
{\it 
$\{ 0 \} 
= \mathcal{O}_{(G(K')^{\circ})_{E_1}}(0)$. 
}

\medskip

(2) 
{\it 
$\mathcal{J}_3(K^{\mathbb{C}})_{-L^{\times}_{2E_1}} 
= (\cup_{c \in \mathbb{C}} 
\mathcal{S}_2(c, K^{\mathbb{C}})) 
\cup \{ 0 \}$ 
such that 
}

\smallskip

\noindent 
(\roman{f1}) 
{\it 
$\mathcal{S}_2(c, K^{\mathbb{C}}) = 
\mathcal{O}_{(G(K^{\mathbb{C}})^{\circ})_{E_1}}
(\sqrt{\frac{c}{2}} (E_2 - E_3))$ 
for $c \in {\mathbb{C}} \backslash \{ 0 \}$; 
}

\smallskip

\noindent 
(\roman{f2}) 
{\it 
$\mathcal{S}_2(0, K^{\mathbb{C}}) 
= \mathcal{O}_{(G(K^{\mathbb{C}})^{\circ})_{E_1}}(M_1)$; 
and 
}

\smallskip

\noindent 
(\roman{f3}) 
{\it 
$\{ 0 \} = 
\mathcal{O}_{(G(K^{\mathbb{C}})^{\circ})_{E_1}}(0)$. 
}

\medskip

{\it Proof.} 
(1) 
For 
$W \in \mathcal{J}_3(K')_{-L^{\times}_{2 E_1}}$, 
put $c := (W | W) \in \mathbb{R}$. 
By Lemma 2.1 (1) and (2) (\roman{f2}), 
$\alpha W = \mathbb{X}(0, s, -s; u \sqrt{-1} e_4, 0, 0)$ 
for some 
$s \geqq 0, u \geqq 0$ 
and $\alpha \in (G(K')^{\circ})_{E_1}^{\tau}$. 
Then $c = (\alpha W | \alpha W) = 2 (s^2 - u^2)$. 
For $t \in \mathbb{R}$, 
put 
$\mathbb{X}(r; x) 
:= \beta_1 (t; \sqrt{-1} e_4, 1) (\alpha W)$, 
so that 
$r_1 = x_2 = x_3 = 0$, 
$r_2 = - r_3 = \cosh (2 t) (s - u \tanh (2t))$ 
and 
$x_1 = v \sqrt{-1} e_4$ 
with 
$v := \cosh (2t) (u - s \tanh (2t))$.

(\roman{f1}-1) 
If $c > 0$, then 
$s > u \geqq 0$ 
and $| u / s | < 1$, 
so that 
$\tanh (2 t) = u / s$ 
for some $t \in \mathbb{R}$ 
such that $v = 0$ 
and $r_2 = \cosh (2 t) (s^2 - u^2)/s > 0$.  
In this case, 
$\mathbb{X}(r; x) = r_2 (E_2 - E_3)$ 
with 
$c = (W | W) 
= (\mathbb{X}(r; x) | \mathbb{X}(r; x)) 
= 2 (r_2)^2$, 
so that 
$\mathbb{X}(r; x) 
= \sqrt{\frac{c}{2}} (E_2 - E_3) 
\in \mathcal{S}_2( c, K' )$.

(\roman{f1}-2) 
If $c < 0$, then $u > s \geqq 0$ 
and $| s / u | < 1$, so that 
$\tanh ( 2 t ) = s / u$ 
for some $t \in \mathbb{R}$ 
such that $r_2 = 0$ and 
$v = \cosh (2t) (u^2 - s^2)/u > 0$. 
In this case, 
$\mathbb{X}(r; x) = v F_1(\sqrt{-1}e_4)$ 
with 
$c = (W | W) 
= (\mathbb{X}(r; x) | \mathbb{X}(r; x)) 
= - 2 v^2$, 
so that 
$\mathbb{X}(r; x) 
= \sqrt{\frac{- c}{2}} F_1 (\sqrt{-1} e_4) 
\in \mathcal{S}_2(c, K')$.

(\roman{f2}, \roman{f3}) 
If $c = 0$, 
then $s^2 - u^2 = c/2 = 0$, 
so that $s = u \geqq 0$ 
and $r_2 = v = u e^{-2t}$. 
When $u \ne 0$: 
$u > 0$ and $u e^{-2 t} = 1$ 
for some $t \in \mathbb{R}$. 
In this case, 
$\mathbb{X}(r; x) 
= E_2 - E_3 + F_1 (\sqrt{-1}e_4)) 
= M_{1'} \in \mathcal{S}_2(0, K')$. 
When $u = 0$: 
$r_2 = v = u = 0$ 
and $\mathbb{X}(r; x) = 0 \in \{ 0 \}$.

(2) 
For 
$W \in \mathcal{J}_3(K^{\mathbb{C}})_{-L^{\times}_{2 E_1}}$, 
put $c := (W|W) \in \mathbb{C}$. 
By Lemma 2.1 (1) and (2) (\roman{f3}), 
$\alpha W 
= \mathbb{X}( 0, t_2 + s_2 \sqrt{-1}, 
- t_2 - s_2 \sqrt{-1}; u, 0, 0)$ 
for some $t_2, s_2, u \in \mathbb{R}$ 
with $s_2, u \geqq 0$ 
and some $\alpha \in (G(K)_{E_1})^{\circ}  
\subseteqq (G(K^{\mathbb{C}})_{E_1})^{\circ}$.  
Then 
$c = (\alpha W | \alpha W) 
= 2((t_2 + s_2 \sqrt{-1})^2 + u^2)$. 
For $t \in \mathbb{R}$, 
put $\mathbb{X}(r; x) := \beta_1(t; 1, \sqrt{-1}) (\alpha W)$ 
with 
$\beta_1(t; 1, \sqrt{-1}) \in (G(K)_{E_1})^{\circ} 
\subseteqq (G(K^{\mathbb{C}})_{E_1})^{\circ}$, 
so that 
$r_1 = x_2 = x_3 = 0$, 
$r_2 = - r_3 = 
(t_2 + s_2 \sqrt{-1}) \cos (2 t) + u \sin (2 t)$ 
and 
$x_1 = u \cos (2 t) - (t_2 + s_2 \sqrt{-1}) \sin (2 t)$.

(\roman{f1}) 
If $c \ne 0$, 
then 
$(t_2 + (s_2 + u) \sqrt{-1}) (t_2 + (s_2 - u) \sqrt{-1}) 
= c/2 \ne 0$, 
so that 
$e^{\sqrt{-1} 4 t} = 
(t_2 + (s_2 + u) \sqrt{-1}) / (t_2 + (s_2 - u) \sqrt{-1}) 
\ne 0$ for some $t \in {\mathbb{C}}$, and that 
$x_1 
= u (e^{\sqrt{-1} 2 t} + e^{-\sqrt{-1} 2 t}) /2 
- (t_2 + s_2 \sqrt{-1}) 
(e^{\sqrt{-1} 2 t} - e^{-\sqrt{-1} 2 t}) / (2 \sqrt{-1}) 
= \frac{\sqrt{-1}}{2} 
\{ (t_2 + (s_2 - u) \sqrt{-1}) e^{\sqrt{-1} 2 t}   
- (t_2 + (s_2 + u) \sqrt{-1}) e^{-\sqrt{-1} 2 t} \} 
= 0$. 
In this case, 
$\mathbb{X}(r; x) = r_2 (E_2 - E_3)$ 
with   
$c = (\mathbb{X}(r; x) | \mathbb{X}(r; x)) 
= 2 (r_2)^2$, 
so that 
$\mathbb{X}(r; x) 
= \sqrt{\frac{c}{2}} (E_2 - E_3) 
\in \mathcal{S}_2(c, K^{\mathbb{C}})$.

(\roman{f2}, \roman{f3}) 
If $c = 0$, then 
$t_2^2 - s_2^2 + u^2 + 2 t_2 s_2 \sqrt{-1} 
= c/2 = 0$, 
so that $t_2 s_2 = 0$.  
When $s_2 = 0$: $t_2 = u = 0$, 
so that $\mathbb{X}(r; x) = 0 \in \{ 0 \}$. 
When $s_2 \ne 0$: 
$t_2 = 0$, 
$u = s_2 > 0$, 
$r_2 = - r_3 = \sqrt{-1} u e^{-2 t \sqrt{-1}}$, 
$x_1 = u e^{-2 t \sqrt{-1}}$. 
There exists $t \in {\mathbb{C}}$ 
such that $\sqrt{-1} u e^{-2 t \sqrt{-1}} = 1$, 
so that 
$\mathbb{X}(r; x) = E_2 - E_3 + F_1(\sqrt{-1}) 
= M_1 \in \mathcal{S}_2(0, K^{\mathbb{C}})$.~\qed

\bigskip

{\sc Lemma 2.3.}
(1) 
{\it 
If  
$Y = \mathbb{X}(r; x) \in \mathcal{J}_2(\tilde{K})$, 
then 
$\mathrm{tr}(Y) = r_2 + r_3$, 
$\mathrm{det}(E_1 + Y) = r_2 r_3 - N(x_1)$ 
and 
$Y^{\times 2} = \mathrm{det}(E_1 + Y) E_1$. 
}

\smallskip

(2) 
{\it 
For any $X \in \mathcal{J}_3(\tilde{K})_{\sigma}$, 
there exists 
$Y \in \mathcal{J}_2(\tilde{K})$ 
such that 
$X = (X | E_1) E_1 + Y$ 
and that 
$Y = \frac{\mathrm{tr}(Y)}{2}(E_2 + E_3) + W$  
for some 
$W \in \mathcal{J}_3(\tilde{K})_{- L^{\times}_{2E_1}}$ 
such that 
$(W, W) = 
\frac{1}{2}(\mathrm{tr}(Y)^2 - 4 \mathrm{det}(E_1 + Y))$. 
In this case, 
put 
$\Psi_Y(\lambda) := 
\lambda^2 - \mathrm{tr}(Y) \lambda 
+ \mathrm{det}(E_1 + Y) \equiv 
(\lambda - \lambda_2)(\lambda - \lambda_3)$ 
with some $\lambda_2, \lambda_3 \in {\mathbb{C}}$. 
Then 
$\Phi_X(\lambda) 
= (\lambda - (X | E_1)) \Psi_Y (\lambda)$ 
and 
$2(W, W) = (\lambda_2 - \lambda_3)^2$. 
}

\smallskip

(3) 
{\it 
$\mathcal{O}_{(G(\tilde{K})^{\tau})^{\circ}}(X) 
\cap \mathcal{J}_2(\tilde{K}) \ne \emptyset$ 
if 
$X \in \mathcal{J}_3(\tilde{K})$ 
with $X^{\times 2} = 0$. 
}

\medskip

{\it Proof.}
(1) 
By Lemma 1.1, one has the first 
and the second equations. 
And 
$Y^{\times 2} 
= \frac{1}{2} (2 r_2 r_3 - 2 N(x_1)) E_1 
= \mathrm{det}(E_1 + Y) E_1$.

(2)
Take 
$X := \mathbb{X}(r_1, r_2, r_3; x_1, 0, 0) 
\in \mathcal{J}_3(\tilde{K})_{\sigma}$. 
Put 

\[
Y := \mathbb{X}(0, r_2, r_3; x_1, 0, 0),
~W := 
\frac{r_2 - r_3}{2} (E_2 - E_3) + F_1(x_1) 
\in \mathcal{J}_2(\tilde{K}). 
\]

\noindent 
Then 
$X = r_1 E_1 + Y$, 
$Y = \frac{r_2 + r_3}{2} (E_2 + E_3) + W$;  
$\mathrm{tr}(Y) = r_2 + r_3$, 
$\mathrm{det}(E_1 + Y) = r_2 r_3 - N(x_1)$ 
and 
$(W | W) = \frac{(r_2 - r_3)^2}{2} + 2 N(x_1) 
= \frac{1}{2} (\mathrm{tr}(Y)^2 - 4~\mathrm{det}(E_1 + Y))$. 
By Lemma 1.1 (1),  
$\varphi_X(\lambda) 
= (\lambda - r_1) (\lambda - r_2) (\lambda-r_3) 
- (\lambda - r_1) N(x_1)
= (\lambda - r_1) (\lambda^2 - (r_2 + r_3) \lambda 
+ (r_2 r_3 - N(x_1))) 
= (\lambda - r_1) \Psi_Y(\lambda)$. 
Because of  
$\Psi_Y(\lambda) \equiv 
(\lambda - \lambda_2) (\lambda - \lambda_3)$, 
one has that 
$\mathrm{tr}(Y) = \lambda_2 + \lambda_3$, 
$\mathrm{det}(E_1 + Y) = \lambda_2 \lambda_3$, 
so that 
$2(W, W) 
= (\lambda_2 + \lambda_3)^2 - 4 \lambda_2 \lambda_3 
= (\lambda_2 - \lambda_3)^2$.

(3) Take $X \in \mathcal{J}_3(\tilde{K})$ 
with $X^{\times 2} = 0$. 
(\roman{f1}) 
When $\tilde{K} = K'$:  
By Lemma 1.5 (4), 
$\alpha X = \sum_{i = 1}^3 (s_i E_i + F_i(p_i \sqrt{-1} e_4))$ 
for some $p_i \in K \cap K'$, $s_i \in \mathbb{R}$ 
and $\alpha \in (G(K')^{\tau})^{\circ}$, 
so that 
$0 = \alpha (X^{\times 2}) = (\alpha X)^{\times 2} 
= \sum_{i=1}^3(s_{i+1}s_{i+2} + 2 N(p_i))E_i 
+\sum_{i=1}^3 F_i(\overline{p_{i+1} p_{i+2}} 
- s_i p_i \sqrt{-1}e_4)$, 
that is, 
$s_{i+1} s_{i+2} + 2 N(p_i) 
= \overline{p_{i+1}}p_{i+2} = s_i p_i = 0$ 
for all $i \in \{ 1, 2, 3 \}$. 
(Case 1) When $p_i = 0$ for all $i$: 
$0 = s_2 s_3 = s_3 s_1 = s_1 s_2$, 
so that 
$\alpha X = s_i E_i$ for some $i$. 
If $i = 2$ or $3$, then 
$\alpha X \in \mathcal{J}_2(K')$. 
If $i = 1$, then 
$\hat{\beta}_{3} (\alpha X) 
= s_1 E_2 \in \mathcal{J}_2(K')$ 
by 
$\hat{\beta}_3 \in (G(K')^{\tau})^{\circ}$ 
defined in Lemma 1.4 (1) (\roman{f2}). 
(Case 2) When $p_i \ne 0$ for some $i$: 
$p_{i+1} = p_{i+2} = 0$ and $s_i = 0$. 
If $i = 1$, then 
$\alpha X \in \mathcal{J}_2(K')$. 
If $i = 2$, then 
$\hat{\beta}_3 (\alpha X) \in \mathcal{J}_2(K')$. 
If $i = 3$, then 
$\hat{\beta}_2 (\alpha X) \in \mathcal{J}_2(K')$ 
by 
$\hat{\beta}_2 \in (G(K')^{\tau})^{\circ}$ 
defined in Lemma 1.4 (1) (\roman{f2}).

(\roman{f2}) 
When 
$\tilde{K} = K^{\mathbb{C}}$: 
$\alpha X = 
Y + \sqrt{-1} \mathrm{diag}(r_1', r_2', r_3')$ 
for some 
$r_i' \in \mathbb{R}$ ($i = 1, 2, 3$), 
$Y \in \mathcal{J}_3(K)$, 
and $\alpha \in G(K)^{\circ} 
= (G(K^{\mathbb{C}})^{\tau})^{\circ}$ 
by Lemma 1.5 (3). 
Putting $Y = \mathbb{X}(r; x)$, 
$s_i := r_i + \sqrt{-1} r_i' \in \mathbb{C}$, 
one has 
$0 = (\alpha X)^{\times 2} 
= \sum_{i = 1}^3 \{ 
(s_{i+1} s_{i+2} - 2 N(x_i)) E_i 
+ F_i(\overline{x_{i+1} x_{i+2}} - s_i x_i) \}$, 
that is, 
$0 
= r_i' x_i 
= \overline{x_{i+1} x_{i+2}}-r_i x_i 
= s_{i+1} s_{i+2}  - 2 N(x_i)$ 
for all $i \in \{ 1, 2, 3 \}$. 
Then  
(Case 1) $x_i = 0$ for all $i$, 
(Case 2) $x_i \ne 0$, $x_{i+1} = x_{i+2} = 0$ 
for some $i$,   
(Case 3) $x_i \ne 0$, $x_{i+1} \ne 0$, $x_{i+2} = 0$ 
for some $i$; or (Case 4) $x_i \ne 0$ for all $i$. 
In (Case 1), $0 = s_{i+1} s_{i+2}$ 
for all $i$, so that 
$\alpha X = s_i E_i$ for some $i$. 
If $i = 2$ or $3$, then 
$\alpha X \in \mathcal{J}_2(K^{\mathbb{C}})$. 
If $i = 1$, then 
$\hat{\beta}_3 (\alpha X) 
= s_1 E_2 \in \mathcal{J}_2(K^{\mathbb{C}})$ 
by $\hat{\beta}_3 \in (G(K^{\mathbb{C}})^{\tau})^{\circ}$ 
defined in Lemma 1.4 (1) (\roman{f2}). 
In (Case 2), 
$0 = r_i' = r_i$, 
so that 
$\alpha X = s_{i+1} E_{i+1} + s_{i+2} E_{i+2} + F_i(x_i)$ 
and that 
$\hat{\beta}_k (\alpha X) 
\in \mathcal{J}_2(K^{\mathbb{C}})$ 
for some 
$\hat{\beta}_k \in (G(K^{\mathbb{C}})^{\tau})^{\circ}$ 
defined in Lemma 1.4 (1) (\roman{f2}). 
In (Case 3), 
$0 = r_i' = r_i 
= r_{i+1}' = r_{i+1} 
= N(x_i) = N(x_{i+1})$, 
so that $\alpha X = s_{i+2} E_{i+2}$ 
and that 
$\hat{\beta}_k (\alpha X) 
\in \mathcal{J}_2(K^{\mathbb{C}})$ 
for some 
$\hat{\beta}_k 
\in (G(K^{\mathbb{C}})^{\tau})^{\circ}$ 
defined in Lemma 1.4 (1) (\roman{f2}). 
In (Case 4), $r_i' = 0$ for all $i$, 
so that $\alpha X \in \mathcal{J}_3(K)$ 
and that $\alpha_1 (\alpha X)$ is diagonal 
for some $\alpha_1 \in G(K)^{\circ}$ 
by Lemma 1.5 (2). 
Then 
$\beta(\alpha_1(\alpha X))) 
\in \mathcal{J}_2(K^{\mathbb{C}})$ 
for some $\beta \in (G(K^{\mathbb{C}})^{\tau})^{\circ}$ 
by the argument on (Case 1).~\qed 

\bigskip

{\it Proof of Proposition 0.1 (3) 
when $\tilde{K} \ne K$.} 
Take any 
$X \in \mathcal{P}_2 (\tilde{K})$. 
By (3), 
$\alpha_1 X 
\in \mathcal{J}_2(\tilde{K})$ 
for some 
$\alpha_1 \in (G(\tilde{K})^{\tau})^{\circ}$. 
By (1), 
$\mathrm{det}(E_1 + \alpha_1 X) E_1 
= (\alpha_1 X)^{\times 2} = \alpha (X^{\times 2}) 
= 0$, i.e. 
$\mathrm{det}(E_1 + \alpha_1 X) = 0$. 
By (2), 
$\alpha_1 X 
= \frac{\mathrm{tr}(\alpha_1 X)}{2} (E_2 + E_3) + W 
= \frac{1}{2}(E_2 + E_3) + W$  
for some 
$W \in \mathcal{J}_3(\tilde{K})_{-L^{\times}_{2E_1}}$ 
such that 
$(W|W) 
= \frac{1}{2}(\mathrm{tr}(\alpha_1 X)^2 
- 4 \mathrm{det}(E_1 + \alpha_1 X)) 
= \frac{1}{2}$, 
so that $W \in \mathcal{S}_2(1/\sqrt{2}, \tilde{K})$. 
By Lemma 2.2 (1)(2), 
$\alpha_2 W = \frac{1}{2}(E_2 - E_3) 
\in \mathcal{S}_2(1/\sqrt{2}, \tilde{K})$ 
for some 
$\alpha_2 \in (G(\tilde{K})_{E_1})^{\circ}$. 
Then 
$\alpha_2 (\alpha_1 X) 
= \frac{1}{2}(E_2 + E_3) + \frac{1}{2}(E_2 - E_3) 
= E_2$ by Lemma 2.1 (1).  
By $\hat{\beta}_3 
\in (G(\tilde{K})_{E_3}^{\tau})^{\circ}$ 
defined in Lemma 1.4 (1) (\roman{f2}), 
$\hat{\beta}_3(\alpha_2 (\alpha_1 X)) = E_1$, 
where 
$\hat{\beta}_3 \alpha_2 \alpha_1 
\in G(\tilde{K})^{\circ}.$~\qed

\bigskip

{\it Proof of Proposition 0.1 (4) (\roman{f1}).} 
Take any $X \in \mathcal{M}_1(\tilde{K})$ 
defined in Lemma 1.6 (5). 
By (3), 
$0 \ne \alpha X \in \mathcal{J}_2(\tilde{K})$ 
for some $\alpha \in (G(\tilde{K})^{\tau})^{\circ}$ 
with $\mathrm{tr}(\alpha X) = \mathrm{tr}(X) = 0$. 
In this case, by (2), 
$\alpha X 
= \frac{\mathrm{tr}(\alpha X)}{2} (E_2 + E_3) + W 
= W$ for some 
$W \in \mathcal{J}_3(\tilde{K})_{-L^{\times}_{2E_1}}$. 
And 
$(\alpha X | \alpha X) 
= (X | X) 
= (X \circ X | E) 
= - 2 (X \times X | E)  
= - 2 \mathrm{tr}(X^{\times 2}) 
= 0$ by Lemma 1.1 (4). 
Hence, 
$\alpha X \in \mathcal{S}_2(0, \tilde{K})$. 
By Lemma 2.2 (1) (\roman{f2}) or (2) (\roman{f2}), 
there exists 
$\beta \in (G(\tilde{K})^{\circ})_{E_1}$ 
such that $\beta (\alpha X) = M_{1'}$ 
(when $\tilde{K} = K'$) or 
$M_1$ (when $\tilde{K} = K^{\mathbb{C}}$).~\qed

\bigskip

\begin{center}
{\bf 3. 
Theorems 0.2 and 0.3 
in (1) (\roman{f1}, \roman{f2}). 
}
\end{center}

\medskip

Assume that $X \in \mathcal{J}_3(\tilde{K})$ 
admits a characteristic root 
$\lambda_1 \in \mathbb{F}$ 
of multiplicity $1$. 
Then 
$0 \ne \Phi_X'(\lambda_1) 
= \mathrm{tr}(\varphi_X (\lambda_1)^{\times 2})$ 
by Lemma 1.2 (2), so that 

\[
E_{X, \lambda_1} 
:= \frac{1}{\mathrm{tr}(\varphi_X(\lambda_1))} \varphi_X(\lambda_1)^{\times 2} 
\in V_X 
\]

\noindent 
is well-defined.  
Put 
$W_{X, \lambda_1} 
:= X - \lambda_1 E_{X, \lambda_1} 
- \frac{\mathrm{tr}(X) -\lambda_1}{2} 
\varphi_{E_{X, \lambda_1}}(1) \in V_X$.
Then 

\[
X = 
\lambda_1 E_{X, \lambda_1} 
+ \frac{\mathrm{tr}(X)-\lambda_1}{2} 
\varphi_{E_{X, \lambda_1}} (1)  
+ W_{X, \lambda_1}. 
\]

\bigskip

{\sc Lemma 3.1.} 
{\it 
Assume that $X \in \mathcal{J}_3(\tilde{K})$ 
admits a characteristic root 
$\lambda_1 \in \mathbb{F}$ 
of multiplicity $1$. Then:}

\medskip

(1) 
{\it 
$V_X \cap \mathcal{P}_2(\tilde{K}) 
\ni E_{X, \lambda_1} \ne 0$, 
$\varphi_{E_{X, \lambda_1}}(1) \ne 0$, 
$E_{X, \lambda_1}^{\times 2} = 0$, 
$2E_{X, \lambda_1}\times \varphi_{E_{X, \lambda_1}}(1) 
= \varphi_{E_{X, \lambda_1}}(1)$, 
$\varphi_{E_{X, \lambda_1}}(1)^{\times 2} 
= E_{X, \lambda_1}$, 
$2E_{X, \lambda_1} \times W_{X, \lambda_1} 
= -W_{X, \lambda_1}$;  
}

\smallskip

(2) 
$V_X = \mathbb{F} E_{X, \lambda_1} 
\oplus \mathbb{F} \varphi_{E_{X, \lambda_1}} (1) 
\oplus \mathbb{F} W_{X, \lambda_1}$ 
{\it such that} 
$v_X = 2$ 
({\it if} $W_{X, \lambda_1} = 0)$ {\it or} 
$v_X = 3$ 
({\it if} $W_{X, \lambda_1} \ne 0$) 
{\it 
with 
$(E_{X, \lambda_1} | \varphi_{E_{X, \lambda_1}} (1)) 
= (E_{X, \lambda_1} | W_{X, \lambda_1}) 
= (\varphi_{E_{X, \lambda_1}}(1) | W_{X, \lambda_1}) 
= 0$, $(E_{X, \lambda_1} | E_{X, \lambda_1}) = 1$, \\
$(\varphi_{E_{X, \lambda_1}} (1) 
| \varphi_{E_{X, \lambda_1}} (1)) = 2$ 
and 
$(W_{X, \lambda_1} | W_{X, \lambda_1}) 
= \Delta_X (\lambda_1)$. 
}

\medskip

{\it Proof.} 
(1) 
Put 
$Z := \varphi_X (\lambda_1)$ 
and $Y := Z^{\times 2}$, 
so that $Y^{\times 2} = 0$ by Lemma 1.6 (4). 
Then 
$E_{X, \lambda_1} 
= \frac{1}{\mathrm{tr}(Y)} Y$,  
$\mathrm{tr}(E_{X, \lambda_1}) = 1$ 
and $E_{X, \lambda_1}^{\times 2} = 0$, 
so that  
$E_{X, \lambda_1} \in \mathcal{P}_2(\tilde{K}) \cap V_X$. 
Note that 
$\mathrm{tr}(\varphi_{E_{X, \lambda_1}}(1)) 
= \mathrm{tr} (E) - \mathrm{tr}(E_{X, \lambda_1}) 
= 3 - 1 = 2 \neq 0$, 
so that $\varphi_{E_{X, \lambda_1}} (1) \neq 0$. 
By Lemma 1.1 (3), 
$2 E_{X, \lambda_1} \times \varphi_{E_{X, \lambda_1}} (1) 
= 2 E_{X, \lambda_1} \times (E - E_{X, \lambda_1}) 
= 2 E_{X, \lambda_1} \times E 
= \mathrm{tr}(E_{X, \lambda_1}) E - E_{X, \lambda_1} 
= \varphi_{E_{X, \lambda_1}} (1)$ 
and 
$\varphi_{E_{X, \lambda_1}} (1)^{\times 2} 
= (E - E_{X, \lambda_1})^{\times 2} 
= E^{\times 2} - 2 E \times E_{X, \lambda_1} 
= E - \varphi_{E_{X, \lambda_1}} (1) 
= E_{X, \lambda_1}$.  
By direct compuations, 
$W_{X, \lambda_1} = 
\frac{\mathrm{tr}(Z)}{2} \varphi_{E_{X, \lambda_1}} (1) 
- Z$. 
By Lemma 1.6 (2) (\roman{f3}) 
and $\mathrm{det}(Z) = 0$, 
$2 E_{X, \lambda_1} \times Z 
= \frac{2}{\mathrm{tr}(Y)} Z^{\times 2} \times Z 
= \frac{2}{\mathrm{tr}(Y)} 
\frac{-1}{2}( 
\mathrm{tr}(Z) Y + \mathrm{tr}(Y) Z 
- \mathrm{tr}(Y) \mathrm{tr}(Z) E 
+ \mathrm{det}(Z) E ) 
= - \mathrm{tr}(Z) E_{X, \lambda_1} - Z + \mathrm{tr}(Z) E 
= - Z + \mathrm{tr}(Z) \varphi_{E_{X, \lambda_1}} (1)$. 
Hence, 
$2 E_{X, \lambda_1} \times W_{\lambda_1} 
= \frac{\mathrm{tr}(Z)}{2} \varphi_{E_{X, \lambda_1}}(1) 
+ Z - \mathrm{tr}(Z) \varphi_{E_{X, \lambda_1}}(1)  
= Z - \frac{\mathrm{tr}(Z)}{2} \varphi_{E_{X, \lambda_1}}(1) 
=- W_{X, \lambda_1}$.

(2) 
Since $V_X$ is spanned by $E, X, X^{\times 2}$, 
$v_X := \mathrm{dim}_{\mathbb{F}} V_X \leqq 3$. 
If 
$W_{X, \lambda_1} \neq 0$, 
then 
$E_{X, \lambda_1}$, 
$\varphi_{E_{X, \lambda_1}} (1)$, 
$W_{X, \lambda_1}$ 
are eigen-vectors of 
$L^{\times}_{2 E_{X,\lambda_1}}$ 
with different eigen-values $0, 1, -1$, 
i.e. $v_X = 3$. 
If $W_{X, \lambda_1} = 0$, 
then 
$X = 
\lambda_1 E_{X, \lambda_1} 
+ \frac{\mathrm{tr}(X)-\lambda_1}{2} 
\varphi_{E_{X, \lambda_1}} (1)$ 
and 
$X^{\times 2} 
= \frac{\lambda_1(\mathrm{tr}(X)-\lambda_1)}{2} 
\varphi_{E_{X, \lambda_1}}(1)  
+ (\frac{\mathrm{tr}(X)-\lambda_1}{2})^2 
E_{X, \lambda_1}$, 
so that $V_X$ is spanned by 
$E_{X, \lambda_1}$ and $\varphi_{E_{X, \lambda_1}} (1) 
= E - E_{X, \lambda_1}$, 
i.e. $v_X = 2$. 
By Lemmas 1.1 (2) and 1.6 (3), 
$L^{\times}_{2 E_{X, \lambda_1}}$ 
is a symmetric $\mathbb{F}$-linear 
transformation on $(V_X, (* | *))$, 
so that 
$E_{X, \lambda_1}, 
\varphi_{E_{X, \lambda_1}}(1), 
W_{X, \lambda_1}$ 
are orthogonal as 
zero or eigen-vectors of 
$L^{\times}_{2 E_{X, \lambda_1}}$ 
with the different eigen-values. 
By (1) and Lemma 1.1 (4), 
$0 = 2 \mathrm{tr}(E_{X, \lambda_1}^{\times 2}) 
= \mathrm{tr}(E_{X, \lambda_1})^2 
- (E_{X, \lambda_1} | E_{X, \lambda_1}) 
= 1 - (E_{X, \lambda_1} | E_{X, \lambda_1})$, 
so that 
$(E_{X, \lambda_1} | E_{X, \lambda_1}) = 1$ 
and 
$(\varphi_{E_{X, \lambda_1}}(1) 
| \varphi_{E_{X, \lambda_1}}(1)) 
= (E | E) - 2 (E | E_{X, \lambda_1}) 
+ (E_{X, \lambda_1} | E_{X, \lambda_1}) 
= 3 - 2 \mathrm{tr}(E_{X, \lambda_1}) + 1 
= 2$. 
Because of the orthogonality in (1), 
$(X | E_{X, \lambda_1}) = \lambda_1$,  
$(X | \varphi_{E_{X, \lambda_1}}(1)) 
= \mathrm{tr}(X) - \lambda_1$ 
and 
$(W_{X, \lambda_1} | W_{X, \lambda_1}) 
= (X | X) 
+ \lambda_1^2 + 
\frac{(\mathrm{tr}(X) - \lambda_1)^2}{2} 
- 2 \lambda_1 (X | E_{X, \lambda_1}) 
- (\mathrm{tr}(X) - \lambda_1) 
(X | \varphi_{E_{X, \lambda_1}}(1)) 
= (X | X) 
- \frac{3}{2}\lambda_1^2 
+ \mathrm{tr}(X) \lambda_1 
- \frac{1}{2} \mathrm{tr}(X)^2 
= \Delta_X (\lambda_1)$. 
\qed

\bigskip

Note that 
$\Delta_X(\lambda_1) 
= - \frac{1}{2} \{ 
3 \lambda_1^2 - 2 \mathrm{tr}(X) \lambda_1 
+ \mathrm{tr}(X)^2 
- 2 (X | X) \} \in \mathbb{F}$ 
is an invariant on ${\cal O}_{G(\tilde{K})}(X)$ 
if $\lambda_1 \in \mathbb{F}$ 
is a characteristic root of multiplicity $1$ 
for $X \in \mathcal{J}_3(\tilde{K})$.   

\bigskip

{\sc Lemma 3.2.} 
{\it 
Assume that 
$X \in \mathcal{J}_3(\tilde{K})$ 
admits an eigen-value 
$\lambda_1 \in \mathbb{F}$ 
of multiplicity $1$. 
Put 
$\Phi_X (\lambda) 
\equiv \Pi_{i = 1}^3 (\lambda - \lambda_i)$ 
for some 
$\lambda_2, \lambda_3 \in {\mathbb{C}}$ 
with $\lambda_1 \ne \lambda_2, \lambda_3$.  
Then 
${\cal O}_{G(\tilde{K})^{\circ}} 
\ni \lambda_1 E_1 
+ \frac{1}{2}(\mathrm{tr}(X)-\lambda_1) 
(E - E_1) + W$ 
for $W \in \mathcal{J}_3(\tilde{K})_{- L^{\times}_{2E_1}}$ 
such that $(W | W) = \Delta_X (\lambda_1) 
= (\lambda_2 - \lambda_3)^2 / 2$ 
given by $\Lambda_X$ and $v_X$ 
as follows:} 

\medskip

(1) 
{\it 
When $\tilde{K} = K'$ with $\mathbb{F} = \mathbb{R}$: 
}

\smallskip

\noindent 
(\roman{f1}-1) 
{\it 
$W = 
\frac{\sqrt{\Delta_{X}(\lambda_1)}}{\sqrt{2}} 
(E_2 - E_3)$ 
with $\# \Lambda_X = v_X = 3$ 
if $\Delta_X(\lambda_1) > 0$; 
}

\smallskip

\noindent 
(\roman{f1}-2) 
{\it 
$W = \frac{\sqrt{-\Delta_{X}(\lambda_1)}}
{\sqrt{2}} F_1(\sqrt{-1}e_4)$ 
with $\# \Lambda_X = v_X = 3$ if $\Delta_X(\lambda_1) < 0$; 
}

\smallskip

\noindent 
(\roman{f2}) 
{\it 
$W = M_{1'}$ if $\Delta_X (\lambda_1) = 0$ 
with $v_X = 3$; 
}

\smallskip

\noindent 
(\roman{f3}) 
{\it 
$W = 0$ 
if 
$\Delta_X (\lambda_1) = 0$ 
with $v_X = 2$. 
}

\medskip

(2) 
{\it 
When $\tilde{K} = K^{\mathbb{C}}$ 
with $\mathbb{F} = \mathbb{C}$: 
}

\smallskip 

\noindent 
(\roman{f1}) 
{\it 
$W = \frac{w}{\sqrt{2}} 
(E_2 - E_3)$ 
for any $w \in \mathbb{C}$ 
such that $w^2 = \Delta_X(\lambda_1)$ 
with $\# \Lambda_X = v_X = 3$ 
if $0 \ne \Delta_X (\lambda_1) \in {\mathbb{C}}$;  
} 

\smallskip

\noindent 
(\roman{f2}) 
{\it 
$W = M_1$ if $\Delta_X (\lambda_1) = 0$ 
with $v_X = 3$; 
} 

\smallskip 

\noindent 
(\roman{f3}) 
{\it 
$W = 0$ if 
$\Delta_X (\lambda_1) = 0$ 
with $v_X = 2$. 
}

\medskip

{\it Proof.} 
By Lemma 3.1 (1) 
and Proposition 0.1 (3), 
$\alpha E_{X, \lambda_1} = E_1$ 
for some $\alpha \in G(\tilde{K})^{\circ}$, 
so that $\alpha \varphi_{E_{X, \lambda_1}}(1) 
= \alpha (E - E_{X, \lambda_1}) 
= E - E_1$. 
Put 
$W' := \alpha W_{X, \lambda_1}$. 
Then 
$\alpha X 
= \lambda_1 E_1 
+ \frac{\mathrm{tr}(X)-\lambda_1}{2}  
(E - E_1) + W'$ 
with 
$\Phi_{\alpha X}(\lambda) 
= \Pi_{i = 1}^3 (\lambda - \lambda_i)$. 
And 
$2 E_1 \times W' 
= \alpha(2E_{X, \lambda_1} \times W_{X, \lambda_1}) 
= -\alpha W_{X, \lambda_1} = - W'$ 
by Lemma 3.1 (1), 
i.e. 
$W' \in 
\mathcal{J}_3(\tilde{K})_{- L^{\times}_{2E_1}} 
\subset \mathcal{J}_3(\tilde{K})_{\sigma}$. 
By Lemmas 3.1 (2) and 2.3 (2), 
$\Delta_X (\lambda_1) 
= (W_{X, \lambda_1} | W_{X, \lambda_1}) 
= (W' | W') = (\lambda_2 - \lambda_3)^2 / 2$, 
which is determined by $\Lambda_X$, 
so that 
$W' \in \mathcal{S}(\Delta_X (\lambda_1)) 
\cup \{ 0 \}$. 
Note that 
$W' = 0$ (or $W' \ne 0$) 
iff $W_{X, \lambda_1} = 0$ 
(resp. $W_{X, \lambda_1} \ne 0$)  
iff $v_X = 2$ (resp. $v_X = 3$) 
by Lemma 3.1 (2). 
If $\Delta_X (\lambda_1) \ne 0$,  
then $(W' | W') \ne 0$, 
so that $W' \ne 0$ and  
$\lambda_2 \ne \lambda_3$, 
i.e. 
$v_X = \# \Lambda_X = 3$: 
By Lemma 2.2 (1) (\roman{f1}-1, 2) 
or (2) (\roman{f1}), 
$W := \beta W'$ 
is given as (1) (\roman{f1}-1, 2) 
or (2) (\roman{f1}) 
for some 
$\beta \in (G(\tilde{K})_{E_1})^{\circ}$, 
so that 
$\beta (\alpha X) 
= 
\lambda_1 \beta E_1 
+ \frac{\mathrm{tr}(X)-\lambda_1}{2}  
\beta(E - E_1) + W 
= 
\lambda_1 E_1 
+ \frac{\mathrm{tr}(X)-\lambda_1}{2}  
(E - E_1) + W$. 
If $\Delta_X (\lambda_1) = 0$, 
then $(W' | W') = 0$, 
so that 
$W' \in \mathcal{S}_2(0, \tilde{K}) \cup \{ 0 \}$: 
By Lemma 2.2 (1) (\roman{f2}, \roman{f3}) 
or (2) (\roman{f2}, \roman{f3}), 
$W := \beta W'$ 
is given as (1) (\roman{f2}, \roman{f3}) 
or (2) (\roman{f2}, \roman{f3}) 
for some 
$\beta \in (G(\tilde{K})_{E_1})^{\circ}$, 
so that 
$\beta (\alpha X) 
= 
\lambda_1 \beta E_1 
+ \frac{\mathrm{tr}(X)-\lambda_1}{2}  
\beta(E - E_1) + W 
= 
\lambda_1 E_1 
+ \frac{\mathrm{tr}(X)-\lambda_1}{2}  
(E - E_1) + W$.~\qed 

\bigskip

{\it Proof of Theorems 0.2 and 0.3 in 
(1) (\roman{f1}, \roman{f2}).} 
Let 
$X \in \mathcal{J}_3(\tilde{K})$ 
be such as $\#\Lambda_X \ne 1$, 
that is, 
$X$ admits no characteristic root 
of multiplicity 3. 
Since the degree of $\Phi_X(\lambda)$ 
equals $3 = 1 + 1 + 1 = 1 + 2$, 
there exists a characteristic root 
$\mu_1 \in \mathbb{C}$ 
of multiplicity $1$. 
If $\mathbb{F} \ni \mu_1$, 
put $\lambda_1 := \mu_1$. 
If $\mathbb{F} \not\ni \mu_1$, 
then $\mathbb{F} = \mathbb{R} \not\ni \mu_1$, 
so that  
$\Phi_X(\lambda) = 
(\lambda - \mu_1) 
(\lambda - \overline{\mu_1}) (\lambda - \nu_1)$ 
for some $\nu_1 \in \mathbb{R}$. 
In this case, put $\lambda_1 := \nu_1$. 
In all cases, put 
$\Lambda_X = \{ \lambda_1, \lambda_2, \lambda_3 \}$ 
with $\# \Lambda_X = 3$ or $2$ 
such that 
$\Phi_X'(\lambda_1) \ne 0$ 
and 
$\mathrm{tr}(X) = \sum_{i = 1}^3 \lambda_i$, 
so that 
$\Delta_X(\lambda_1) = (\lambda_2 - \lambda_3)^2/ 2$ 
by Lemmas 3.1 (2) and 2.3 (2). 
By Lemma 3.2 (2) 
(if $\tilde{K} = K^{\mathbb{C}}$) 
or (1) (if $\tilde{K} = K'$), 
$\alpha X = 
\lambda_1 E_1 
+ \frac{\mathrm{tr}(X) - \lambda_1}{2} (E - E_1) 
+ W$ 
for some 
$W \in \mathcal{J}_3(\tilde{K})_{-L^{\times}_{2 E_1}}$ 
and $\alpha \in G(\tilde{K})^{\circ}$.

\medskip

(0.2.1) 
{\it When $\mathbb{F} = \mathbb{C}$:} 
$\tilde{K} = K^{\mathbb{C}} 
= \mathbb{R}^{\mathbb{C}}, 
\bC^{\mathbb{C}}, \bH^{\mathbb{C}}$ 
or 
$\bO^{\mathbb{C}}$. 

\smallskip

(0.2.1.\roman{f1}) 
{\it The case of 
$\# \Lambda_X = 3$:} 
Put 
$w := (\lambda_2 - \lambda_3)/ \sqrt{2}$. 
Then 
$\lambda_2 \ne \lambda_3$. 
And 
$\Delta_X(\lambda_1) = w^2 \ne 0$. 
By Lemma 3.2 (2) (\roman{f1}), 
$v_X = 3$ 
and 
$\alpha X = 
\lambda_1 E_1 
+ \frac{\mathrm{tr}(X) - \lambda_1}{2}(E - E_1) 
+ \frac{w}{\sqrt{2}}(E_2 - E_3) 
= \lambda_1 E_1 
+ \frac{\lambda_2 + \lambda_3}{2}(E_2 + E_3) 
+ \frac{\lambda_2 - \lambda_3}{2}(E_2 - E_3) 
= \mathrm{diag}(\lambda_1, \lambda_2, \lambda_3)$. 

\smallskip

(0.2.1.\roman{f2}) 
{\it The case of 
$\# \Lambda_X = 2$:} 
$\lambda_2 = \lambda_3$ 
and $\Delta_X(\lambda_1) = 0$.

(0.2.1.\roman{f2}-1) 
{\it When $v_X = 2$:} 
By Lemma 3.2 (2) (\roman{f3}), 
$\alpha X = \lambda_1 E_1 + \lambda_2(E_2 + E_3) 
= \mathrm{diag} (\lambda_1, \lambda_2, \lambda_2)$.

(0.2.1.\roman{f2}-2) 
{\it When $v_X = 3$:} 
By Lemma 3.2 (2) (\roman{f2}), 
$\alpha X = \lambda_1 E_1 + \lambda_2(E_2 + E_3) + M_1 = 
\mathrm{diag} (\lambda_1, \lambda_2, \lambda_2) + M_1$.

\medskip

(0.3.1) 
{\it When $\mathbb{F} = \mathbb{R}$:} 
$\tilde{K} = K' 
= \bC', \bH'$ 
or $\bO'$. 
And 
$\lambda_1 \in \mathbb{R}$, 
$\lambda_2, \lambda_3 \in \mathbb{C}$. 

\smallskip

(0.3.1.\roman{f1}) 
{\it The case of $\# \Lambda_X = 3$:} 
 
(0.3.1.\roman{f1}-1) 
{\it When $\Lambda_X \subset \mathbb{R}$:} 
It can be assumed that 
$\lambda_1 > \lambda_2 > \lambda_3$ 
by translation if necessary. 
Then 
$\Delta_X(\lambda_1) > 0$. 
By Lemma 3.2 (1) (\roman{f1}-1), 
$v_X = 3$ and 
$\alpha X = \lambda_1 E_1 
+ \frac{\lambda_2 + \lambda_3}{2}(E_2 + E_3) 
+ \frac{\lambda_2 - \lambda_3}{2}(E_2 - E_3)
= \mathrm{diag}(\lambda_1, \lambda_2, \lambda_3)$. 

(0.3.1.\roman{f1}-2) 
{\it When $\Lambda_X \not\subset \mathbb{R}$:}  
$\{ \lambda_2, \lambda_3 \} 
= \{ p \pm q \sqrt{-1} \}$ 
for some $p, q \in \mathbb{R}$ with $q > 0$. 
And $\Delta_{X}(\lambda_1) = -2 q^2 < 0$. 
By Lemma 3.2 (1) (\roman{f1}-2), 
$\alpha X 
= \lambda_1 E_1 + p(E_2 + E_3) 
+ q F_1(\sqrt{-1}e_4) 
= \mathrm{diag}(\lambda_1, p, p) 
+ F_1(q \sqrt{-1}e_4)$. 

\smallskip

(0.3.1.\roman{f2}) 
{\it The case of $\# \Lambda_X = 2$:} 
$\Lambda_X = \{ \lambda_1, \lambda_2 \}$ 
with $\Phi_X'(\lambda_2) = 0$. 
Then 
$\lambda_2 
= \frac{1}{2} (\mathrm{tr}(X) - \lambda_1) 
\in \mathbb{R}$ 
and 
$\Delta_X(\lambda_1) = 0$.

(0.3.1.\roman{f2}-1) 
{\it When $v_X = 2$:} 
By Lemma 3.2 (1) (\roman{f3}), 
$\alpha X 
= \lambda_1 E_1 + \lambda_2 (E_2 + E_3) 
= \mathrm{diag}(\lambda_1, \lambda_2, \lambda_2)$.

(0.3.1.\roman{f2}-2) 
{\it When $v_X = 3$:} 
By Lemma 3.2 (1) (\roman{f2}), 
$\alpha X = \lambda_1 E_1 + \lambda_2(E_2 + E_3) + M_{1'} 
= \mathrm{diag}(\lambda_1, \lambda_2, \lambda_2) + M_{1'}$.~\qed

\bigskip

\begin{center}
{\bf 4. 
Proposition 0.1 (4) 
and Theorems 0.2 and 0.3 in (1) (\roman{f3}).} 
\end{center} 
 
\medskip

Assume that $\tilde{K} \ne K$, i.e., 
$\tilde{K} = \bR^{\mathbb{C}}, 
\bC^{\mathbb{C}}, \bH^{\mathbb{C}}, 
\bO^{\mathbb{C}}; \bC', \bH'$ or $\bO'$. 
Put 
$\mathcal{N}_1(\tilde{K}) 
:= (\mathcal{J}_3(\tilde{K})_0)_{L_{\tilde{M}_1}, 0}$ 
and 
$\mathcal{N}_2(\tilde{K}) 
:= \{ X \in \mathcal{J}_3(\tilde{K})_0 
|~X^{\times 2} = \tilde{M}_1 \}$. 

\bigskip

{\sc Lemma 4.1.}
(1) 
$\mathcal{N}_2(\tilde{K}) 
\subseteqq \mathcal{N}_1(\tilde{K})$. 

\smallskip

(2) 
$\mathcal{N}_1(\tilde{K}) 
= \{ \tilde{M}_1(x) + \tilde{M}_{2 3}(y) 
|~x, y \in \tilde{K} \}$;   

\smallskip

(3) 
$\mathcal{N}_2(\tilde{K})  
= \{ s \tilde{M}_1 + \tilde{M}_{2 3}(y) 
|~s \in \mathbb{F}, 
y \in \mathcal{S}_1(1, \tilde{K}) \}$;  

\medskip

{\it Proof.} 
(1) 
Take $X \in \mathcal{N}_2(\tilde{K})$. 
Then $\mathrm{tr}(X) = 0$ 
and $\mathrm{tr}(X^{\times 2}) 
= \mathrm{tr}(\tilde{M}_1) = 0$. 
By Lemma 1.6 (2) (\roman{f2}), 
$\mathrm{det}(X) X 
= (X^{\times 2})^{\times 2} 
= \tilde{M}_1^{\times 2} = 0$, 
so that $\mathrm{det}(X) = 0$. 
By Lemma 1.6 (2) (\roman{f3}), 
$X \times \tilde{M}_1 
= X \times (X^{\times 2}) 
= - \frac{1}{2}(\mathrm{tr}(X) X^{\times 2} 
+ \mathrm{tr}(X^{\times 2}) X
- (\mathrm{tr}(X) \mathrm{tr}(X^{\times 2}) 
- \mathrm{det}(X))E) 
= 0$, 
as required.

(2) 
(\roman{f1}) 
Take $X = \mathbb{X}(r; x) 
\in \mathcal{N}_1(K^{{\mathbb{C}}})$. 
Then $r_1 + r_2 + r_3 = 0$ 
and that 
$0 = 2 M_1 \times X 
= 2 \mathbb{X}(0, 1, -1; \sqrt{-1}, 0, 0) 
\times \mathbb{X}(r; x) 
= (r_3 - r_2 - 2 (\sqrt{-1} | x_1)) E_1 
- r_1 E_2 + r_1 E_3 
+ F_1(-\sqrt{-1} r_1)
+ F_2(\sqrt{-1} \overline{x_3} - x_2) 
+ F_3(\sqrt{-1} \overline{x_2} + x_3)$ 
by Lemma 1.1 (1), that is, 
$x_2 = \sqrt{-1} \overline{x_3}$, 
$r_1 = 0$, 
$r_2 = - r_3 = -(\sqrt{-1} | x_1) 
= (1 | -\sqrt{-1} x_1)$ 
with 
$x_1 = \sqrt{-1} (-\sqrt{-1} x_1)$, 
so that 
$X = M_1(-\sqrt{-1} x_1) + M_{2 3}(x_3)$.

(\roman{f2}) 
Take $X = \mathbb{X}(r; x) 
\in \mathcal{N}_2(K')$. 
Then $r_1 + r_2 + r_3 = 0$ 
and that 
$0 = 2 M_{1'} \times X 
= 2 \mathbb{X}(0, 1, -1; \sqrt{-1}e_4, 0, 0) 
\times \mathbb{X}(r; x) = 
(r_3 - r_2 - 2 (\sqrt{-1}e_4 | x_1)) E_1 
- r_1 E_2 + r_1 E_3 
+ F_1 (- r_1 \sqrt{-1} e_4)
+ F_2 (- \sqrt{-1} e_4 \overline{x_3} - x_2) 
+ F_3(- \overline{x_2} \sqrt{-1} e_4 + x_3)$ 
by Lemma 1.1 (1), that is, 
$x_2 = - \sqrt{-1} e_4 \overline{x_3}$, 
$r_1 = 0$, 
$r_2 = - r_3 = - (\sqrt{-1} e_4 | x_1) 
= (1 | \sqrt{-1} e_4 x_1)$ 
with 
$x_1 = (\sqrt{-1} e_4)^2 x_1 
= \sqrt{-1} e_4 (\sqrt{-1} e_4 x_1)$, 
so that 
$X = M_{1'}(\sqrt{-1} e_4 x_1) + M_{2' 3}(x_3)$.

(3) (\roman{f1}) 
Take $X \in \mathcal{N}_2(K^{\mathbb{C}})$. 
By (1), $X \in \mathcal{N}_1(K^{\mathbb{C}})$. 
By (2), 
$X = M_1(x_1) + M_{2 3}(x_3)$ 
for some $x_1, x_3 \in K^{\mathbb{C}}$, so that 
$X^{\times 2} = 
M_1(x_1)^{\times 2} 
+ 2 M_1(x_1) \times M_{2 3}(x_3) 
+ M_{2 3}(x_3)^{\times 2} 
= \sqrt{-1} \{ (x_1|1)^2 - N(x_1) \} E_1 
- F_2 (\sqrt{-1} (\overline{x_1} - (x_1|1)) \overline{x_3}) 
- F_3 (x_3 (\overline{x_1} - (x_1 | 1))) 
+ N(x_3) M_1$. 
Hence, 
$X^{\times 2} = M_1$ 
iff 
$N (x_3) = 1$ 
and $\overline{x_1} = (x_1|1) 
\in \mathbb{R}^{\mathbb{C}}$, 
i.e. 
$(x_1, x_3) = (s, x)$ 
for some 
$(s, x) \in 
\mathbb{R}^{\mathbb{C}} \times 
\mathcal{S}_1(1, K^{\mathbb{C}})$.

(\roman{f2}) 
Take $X \in \mathcal{N}_2(K')$. 
By (1), $X \in \mathcal{N}_1(K')$. 
By (2), 
$X = M_{1'}(x_1) + M_{2' 3}(x_3)$ 
for some $x_1, x_3 \in K'$, 
so that 
$X^{\times 2} = 
M_{1'}(x_1)^{\times 2} 
+ 2 M_{1'}(x_1) \times M_{2' 3}(x_3) 
+ M_{2' 3}(x_3)^{\times 2} 
= (N(x_1) - (x_1|1)^2) E_1 
- F_2(((\overline{x_1} 
- (x_1|1))\sqrt{-1}e_4) \overline{x_3}) 
+ F_3(-(x_3 \sqrt{-1}e_4)(\overline{x_1}\sqrt{-1}e_4) 
+(x_1|1)x_3) + N(x_3) M_{1'}$. 
Hence, 
$X^{\times 2} = M_{1'}$ 
iff $N(x_3) = 1$ and $x_1 \in \bR$, as required.~\qed

\bigskip

{\sc Lemma 4.2.} 
$\mathcal{N}_2(\tilde{K}) 
= \mathcal{O}_{(G(\tilde{K})^{\circ})_{\tilde{M}_1}}
(\tilde{M}_{2 3})$.

\medskip

{\it Proof}. 
(1) 
Take any $X \in \mathcal{N}_2(K^{\mathbb{C}})$. 
By Lemma 4.1 (3), 
$X = s M_1 + M_{2 3}(x)$ 
for some $s \in \mathbb{R}^{\mathbb{C}}$ 
and $x \in \mathcal{S}_1(1, K^{\mathbb{C}})$. 
Put  
$x = \sum_{i = 0}^{d_K -1} \xi_i e_i$ 
with $\xi_i \in \mathbb{R}^{\mathbb{C}}$ 
and $x_v := \sum_{i = 1}^{d_K -1} \xi_i e_i$ 
such that $\overline{x_v} = - x_v$.

(Case 1) 
When $\xi_0 = 0$: 
By $x_v = x \in \mathcal{S}_1(1, K^{\mathbb{C}})$, 
$x x = - x \overline{x} = - 1$.  
By Lemma 1.4 (2) (\roman{f2}), 
$\delta_1(x) X 
= s (E_2 - E_3 - F_1 (\sqrt{-1}))  
- F_2(\sqrt{-1}) + F_3(1)$, 
so that 
$\sigma_3 \delta_1(x) X = s M_1 + M_{2 3}$. 
By Lemma 1.4 (1) (\roman{f3}), 
$\beta_{2 3}(t) (\sigma_3 \delta_1(x) X) 
= (2 t + s) M_1 + M_{2 3} = M_{2 3}$ 
if $t = -s/2$ with 
$\beta_{2 3}(-s/2) 
\in (G(K^{\mathbb{C}})^{\circ})_{M_1}$ 
and 
$(\beta_{2 3}(-s/2) \sigma_3 \delta_1(x)) M_1 
= \beta_{2 3}(-s/2) (M_1) = M_1$, 
as required.

(Case 2) 
When $\xi_0 \ne 0$: 
Put 
$Y := \beta_{2 3}(-s/(2 \xi_0)) X$ 
with 
$\beta_{2 3}(-s/(2 \xi_0)) 
\in (G(K^{\mathbb{C}})^{\circ})_{M_1}$. 
By Lemma 1.4 (1) (\roman{f3}), 
$Y = M_{2 3}(x)$.

(\roman{f1}) 
When $x_v = 0$: 
$\xi_0^2 = (x|x) - (x_v|x_v) = 1 - 0 = 1$, 
$x = \xi_0 = \pm 1$ and 
$Y = M_{2 3}(\pm 1) = \pm M_{2 3}$, 
so that 
$Y = M_{2 3}$ or  
$\sigma_1 Y = M_{2 3}$ 
with 
$\sigma_1 \in (G(K^{\mathbb{C}})^{\circ})_{M_1}$.

(\roman{f2}) 
When $x_v \ne 0$ with $d_K \geqq 4$: 
Then $d_K-1 \geqq 3$. 
If $\xi_1^2 + \xi_j^2 = 0$  
for all $j \in \{ 2, \cdots, d_K - 1 \}$ 
and $\xi_2^2 + \xi_3^2 = 0$, then 
$- \xi_1^2 = \xi_2^2 = \xi_3^2 = \cdots = \xi_{d_K-1}^2 = 0$, 
that is, $x_v = 0$, a contradiction. 
Hence, $\xi_i^2 + \xi_j^2 \ne 0$ 
for some $i, j \in \{ 1, \cdots, d_K - 1 \}$ 
with $i \ne j$. 
Take $c \in \mathbb{R}^{\mathbb{C}}$ 
such that $c^2 = \xi_i^2 + \xi_j^2$.  
Put $a := (\xi_j e_i - \xi_i e_j)/c$, 
so that 
$- \overline{a} = a 
\in \mathcal{S}_1(1, K^{\mathbb{C}})$ 
and $(a | x) = 0$. 
Put $y := x \overline{a}$. 
Then $(y | 1) = (a | x) = 0$ and 
$\sigma_3 \delta_1^a Y = F_2(\sqrt{-1} \overline{y}) + F_3(y)$ 
with 
$\sigma_3 \delta_1^a \in (G(K^{\mathbb{C}})^{\circ})_{M_1}$. 
According to the (Case 1), 
$\beta (\sigma_3 \delta_1^a Y) = M_{2 3}$ 
for some $\beta \in (G(K^{\mathbb{C}})^{\circ})_{M_1}$.

(\roman{f3}) 
When $x_v \ne 0$ with $d_K \leqq 4$: 
By Lemma 1.4 (2) (\roman{f3}), 
$\beta_1(x) Y = M_{2 3}$ 
with $\beta_1(x) \in (G(K^{\mathbb{C}})^{\circ})_{M_1}$.

(2) 
Take any $X \in \mathcal{N}_2(K')$. 
Then 
$X = s M_{1'} + M_{2' 3}(x)$ 
for some $s \in \mathbb{R}$ 
and $x \in \mathcal{S}_1(1, K')$ 
by Lemma 4.1 (3). 
Take $p, q \in K_{\tau}'$ 
such that 
$x = p + q \sqrt{-1} e_4$, 
so that 
$\overline{p} p - \overline{q} q = N(x) = 1$. 
Note that 
$\sqrt{-1} e_4 x = 
\overline{q} + \overline{p} \sqrt{-1} e_4$, 
and that   
$x (\sqrt{-1} e_4 x) 
= (p + q \sqrt{-1} e_4) 
(\overline{q} + \overline{p} \sqrt{-1} e_4) 
= p(q + \overline{q}) + (q^2 + \overline{p} p) 
\sqrt{-1} e_4$.

(Case 1) 
When $(\sqrt{-1} e_4 x | 1) = 0$: 
Then $q + \overline{q} = 2(\overline{q} | 1) = 0$. 
By $q = - \overline{q}$, 
$q^2 + \overline{p} p = 
\overline{p} p - \overline{q} q = 1$, 
so that $x (\sqrt{-1} e_4 x) = \sqrt{-1} e_4$ 
and 
$\overline{x} \sqrt{-1} e_4 \overline{x} 
= - \overline{x \sqrt{-1} e_4 x} 
= - \overline{\sqrt{-1} e_4} 
= \sqrt{-1} e_4$.     
By Lemma 1.4 (2) (\roman{f2}), 
$\delta_1(x) M_{1'} 
= E_2 - E_3 + F_1 (x \sqrt{-1} e_4 x) 
= M_{1'}$ 
and 
$\delta_1(x) X 
= s M_{1'} 
+ F_2(-\overline{x} \sqrt{-1} e_4 \overline{x}) 
+ F_3(x \overline{x}) 
= s M_{1'} + M_{2' 3}$. 
By Lemma 1.4 (1) (\roman{f3}), 
$\beta_{2' 3}(-s/2) (s M_{1'} + M_{2' 3}) 
= M_{2' 3}$ 
with 
$\beta_{2 3'} (t) \in (G(K')^{\circ})_{M_{1'}}$.

(Case 2) 
When $(\sqrt{-1} e_4 x | 1) \ne 0$: 
Then $(\sqrt{-1} e_4 | x) = - (1 | \sqrt{-1} e_4 x) \ne 0$. 
By Lemma 1.4 (1) (\roman{f3}), 
$\beta_{2 3'}(-s/(2 (\sqrt{-1} e_4 | x))) X 
= M_{2' 3}(x)$ 
with 
$\beta_{2' 3} (t) \in (G(K')^{\circ})_{M_{1'}}$. 
Note that 
$q + \overline{q} = 2 (\sqrt{-1} e_4 |x) \ne 0$, 
so that $q \in K_{\tau}'$ and $q \ne 0$.

(\roman{f1}) 
When $d_{K'} \geqq 4$: 
By ${\rm dim}_{\bR} K_{\tau}' = d_{K'}/2 \geqq 2$, 
there exists $q_1 \in K_{\tau}'$ 
such that $(\overline{q} | q_1) = 0$, 
so that $(\sqrt{-1} e_4 (x \overline{q_1}) | 1) 
= - (x \overline{q_1} | \sqrt{-1} e_4) 
= - (\overline{q_1} | \overline{x} \sqrt{-1} e_4) 
= (q_1 | \sqrt{-1} e_4 x) 
= (q_1 | \overline{q}) = 0$. 
Put 
$a := q_1/\sqrt{N(q_1)} 
\in \mathcal{S}_1(1, K')$. 
Because of 
$(\sqrt{-1} e_4 a | 1) = 0$, 
$\delta_1(a) M_{1'} = M_{1'}$ 
as well as (Case 1), so that 
$\mathcal{N}_2(K') \ni \delta_1(a) M_{2' 3}(x) 
= M_{2' 3}(x \overline{a})$ 
with $x \overline{a} \in \mathcal{S}_1(1, K')$ 
such that $(\sqrt{-1} e_4 (x \overline{a}) | 1) = 0$. 
Then 
$\delta_1(a) M_{1'} = M_{1'}$ 
and 
$\delta_1(x \overline{a}) M_{2' 3}(x \overline{a}) 
= M_{2' 3}$ 
as well as (Case 1).

(\roman{f2}) 
When $d_{K'} \leqq 2$: 
Then $K' = \bC'$, 
so that $x \in \mathcal{S}_1(1, \bC')$. 
By Lemma 1.4 (2) (\roman{f3}), 
$\beta_1(x) \in (G(\bC')^{\circ})_{M_{1'}}$ 
such that 
$\beta_1(x) M_{2' 3}(x) = M_{2' 3}$. 
\qed

\bigskip

{\it Proof of Proposition 0.1 (4) (\roman{f2}).} 
Take any 
$X \in \mathcal{M}_{2 3}(\tilde{K})$. 
By Lemma 1.6 (5), 
$X^{\times 2} \in \mathcal{M}_1(\tilde{K})$. 
By Proposition 0.1 (4) (\roman{f1}), 
there exists 
$\alpha \in G(\tilde{K})^{\circ}$ 
such that 
$\tilde{M}_1 
= \alpha (X^{\times 2}) = (\alpha X)^{\times 2}$, 
so that 
$\alpha X \in \mathcal{N}_2(\tilde{K})$.
By Lemma 4.2, there exists 
$\beta \in (G(\tilde{K})^{\circ})_{\tilde{M}_1}$ 
such that 
$\beta (\alpha X) = \tilde{M}_{2 3}$, 
as required.~\qed

\bigskip

{\it Proof of Theorems 0.2 and 0.3 
in (1) (\roman{f3}).} 
Take any 
$X \in \mathcal{J}_3(\tilde{K})$ 
with $\# \Lambda_X = 1$ such as 
$\Lambda_X = \{ \lambda_1 \}$. 
By Lemmas 1.2 (1) and (3), 
$X = \lambda_1 E + X_0$ 
with some 
$X_0 \in \{ 0 \} \cup \mathcal{M}_1(\tilde{K}) 
\cup \mathcal{M}_{2 3}(\tilde{K})$.

(\roman{f3}-1) 
When $X_0 = 0$:  
$X = \lambda E$ and 
$v_X = \mathrm{dim}_{\mathbb{F}} V_X 
= \mathrm{dim}_{\mathbb{F}} 
\{ a X^{\times} + b X + c E |~a, b, c \in \mathbb{F} \} 
= \mathrm{dim}_{\mathbb{F}} \{ c E |~c \in \mathbb{F} \} 
= 1$.

(\roman{f3}-2) 
When $X_0 \in \mathcal{M}_1(\tilde{K})$: 
By Proposition 0.1 (4) (\roman{f1}), 
there exists 
$\alpha \in G(\tilde{K})^{\circ}$ 
such that 
$\alpha X = \lambda_1 E + \tilde{M}_1$. 
By Lemma 1.1 (3) with 
$\tilde{M}_1^{\times 2} 
= \mathrm{tr}(\tilde{M}_1) = 0$, 
one has that 
$(\alpha X)^{\times 2} 
= \lambda_1^2 E - \lambda_1 \tilde{M}_1$, 
so that 
$v_X = \mathrm{dim}_{\mathbb{F}} 
\{ a (\alpha X)^{\times 2} + b \alpha X + c E 
|~a, b, c \in \mathbb{F} \} 
= \mathrm{dim}_{\mathbb{F}} 
\{ (a \lambda_1^2 + b \lambda_1 + c) E 
+ (b - a \lambda_1) \tilde{M}_1 |~a, b, c \in \mathbb{F} \} 
= 2$.

(\roman{f3}-3) 
When $X_0 \in \mathcal{M}_{2 3}(\tilde{K})$: 
By Proposition 0.1 (4) (\roman{f2}), 
there exists 
$\alpha \in G(\tilde{K})^{\circ}$ 
such that 
$\alpha X = \lambda_1 E + \tilde{M}_{2 3}$. 
By Lemma 1.1 (3) with 
$\tilde{M}_{2 3}^{\times 2} = \tilde{M}_1$ 
and $\mathrm{tr}(\tilde{M}_{2 3}) = 0$, 
one has that 
$(\alpha X)^{\times 2} 
= \lambda_1^2 E - \lambda_1 \tilde{M}_{2 3} + \tilde{M}_1$, 
so that 
$v_X = 
\mathrm{dim}_{\mathbb{F}} 
\{ a (\alpha X)^{\times 2} 
+ b \alpha X + c E |~a, b, c \in \mathbb{F} \} 
= \mathrm{dim}_{\mathbb{F}} 
\{ (a \lambda_1^2 + b \lambda_1 + c) E 
+ (b - a \lambda_1) \tilde{M}_{2 3} + a \tilde{M}_1 
|~a, b, c \in \mathbb{F} \} = 3$.~\qed

\bigskip

{\bf Acknowledgements.}~The authors are grateful 
to Prof. Ichiro Yokota for introducing them 
to a comprehensive study program on 
explicit realizations of 
the identity connected component of 
exceptional simple Lie groups 
and their subgroups. 
The authors are also grateful 
to the anonymous referee for 
the helpful comment 
to revise the manuscript. 

\bigskip

\renewcommand{\refname}{References}

\begin{flushright}
Graduate School of Engineering, \\
University of Fukui, \\
Fukui-shi, 910-8507, Japan 
\end{flushright}


\begin{thebibliography}{nishioyasukura}
\bibitem{Cc} 
{\sc C.~Chevalley}, 
\textit{Theory of Lie Groups}, 
Princeton UP, 1946. 


\bibitem{CS1950} 
{\sc C.~Chevalley} and {\sc R.D.~Schafer}, 
The exceptional simple Lie algebras $F_4$ and $E_6$, 
{\it Proc. Nat. Acad. Sci. U.S.} {\bf 36} (1950), 
137--141. 

\bibitem{CS2003} 
{\sc J.H.~Conway} and {\sc D.A.~Smith}, 
{\it On Quaternions and Octanions: Their Geometry, 
Arithmetic and Symmetry}, 
A.K.Peters, Ltd, 2003. 
 
\bibitem{Dle1919}
{\sc L.E.~Dickson}, 
On quaternions and their generalizations and the history 
of the eight square theorem, 
{\it Annals of Math.} {\bf 20} (1919), 155--171, 297. 


\bibitem{Dd1978} 
{\sc D.~Drucker}, 
\textit{Exceptional Lie Algebras and the 
Structures of Hermitian Symmetric Spaces}, 
Memoirs of the AMS {\bf 208}, 1978, AMS. 


\bibitem{FK1994}
{\sc J.~Faraut and A.~Koranyi}, 
\textit{Analysis on Symmetric Cones}, 
Oxford UP, New York, 1994.


\bibitem{Fh1951}
{\sc H. Freudenthal}, 
\textit{Oktaven, Ausnahmergruppen und Oktavengeometrie}, 
Math. Inst. Rijksuniv. te Utrecht, 1951; 
(a reprint with corrections, 1960) 
= \textit{Geometriae Dedicata} {\bf 19}-1 (1985), 
7--63. 

\bibitem{Fh1953} 
{\sc H. Freudenthal}, 
Zur ebenen Oktavengeometrie, 
 \textit{Nedel. Akad. Wetensch. Proc. Ser. A.} {\bf 56} 
= \textit{Indag. Math.} {\bf 15} (1953), 195--200,


\bibitem{Fh1964} 
{\sc H. Freudenthal}, 
Lie groups in the foundations of geometry, 
\textit{Adv. Math.} {\bf 1} (1964), 145--190. 

\bibitem{GOV1994} 
{\sc V.V. Gorbatsevich, A.L. Onishchik, E.B. Vinberg}, 
{\it 
Structure of Lie Groups and Lie Algebras \Roman{f3}}, 
Encyclopedia of Math. Sci. {\bf 41}, Springer-Verlag, 1994. 
(Translated from the Russian edition (1990, Moscow) 
by V. Minachin)


\bibitem{Hfr1990}
{\sc F. Reese Harvey}, 
\textit{Spinors and Calibrations}, 
Academic Press, 1990. 


\bibitem{Hs} 
{\sc S.~Helgason}, 
\textit{Differential Geometry, 
Lie Groups, and Symmmetric Spaces}, 
Academic Press, New York 1978. 


\bibitem{Jn1959} 
{\sc N.~Jacobson}, 
Some groups of transformations 
defined by Jordan algebras, \Roman{f1}, 
\textit{J. Reine Angew. Math.} 
{\bf 201}(1959), 178--195.   

\bibitem{Jn1960} 
{\sc N. Jacobson}, 
\textit{Cayley Planes}, 
Mimeographed notes, 
Yale Univ., New Haven, Conn., 1960. 

\bibitem{Jn1961} 
{\sc N.~Jacobson}, 
Some groups of transformations 
defined by Jordan algebras, \Roman{f3}, 
\textit{J. Reine Angew. Math.} 
{\bf 207}(1961), 61--85.  


\bibitem{Jn1968} 
{\sc N.~Jacobson}, 
\textit{Structure and Representations of Jordan Algebras}, 
AMS Colloquim Pub. Vol {\bf XXXIV}, 1968.  

\bibitem{Ks2002} 
{\sc S.~Krutelevich}, 
On a canonical form of 
a $3 \times 3$ Hermitian matrix over 
the ring of integral split octonions,  
\textit{J.Algebra} {\bf 253} (2002), 
276--295. 

\bibitem{My1954} 
{\sc Y. Matsushima}, 
Some remarks on the exceptional simple Lie group $F_4$, 
{\it Nagoya Math. J.} {\bf 4} (1954), 
83--88. 
 

\bibitem{MY2001} 
{\sc T.~Miyasaka} and {\sc I.~Yokota}, 
Constructive diagonalization of an element $X$ of the Jordan 
algebra $\mathcal{J}$ by the exceptional group $F_4$, 
{\it Bull. Fac. Edu 
\& Human Sci. Yamanashi Univ.} {\bf 2}-2 (2001), 
7--10. 
 
\bibitem{Ms1989} 
{\sc S. Murakami}, 
Exceptional simple Lie groups and related topics 
in recent differential geometry, 
{\it Differential Geometry and Topology}, 
LNM {\bf 1369}(1989), Springer, 183--221. 

\bibitem{Rpk1974} 
{\sc P.K. Rasevskii}, 
A theorem on the connectedness of a subgroup 
of a simply connected Lie group commuting with 
any of its automorphisms, 
{\it Trans. Moscow Math. Soc.} {\bf 30} (1974), 
1--24. 

\bibitem{SY1979}
{\sc O. Shukuzawa} and {\sc I. Yokota}, 
Non-compact simple Lie groups $E_{6(6)}$ of type $E_6$, 
\textit{J. Fac. Sci. Shinshu Univ.} {\bf 14}-1 (1979), 
1--13.  


\bibitem{So2006}
{\sc O. Shukuzawa}, 
Explicit classification of orbits 
in Jordan algebra and Freudenthal 
vector space over the exceptional Lie groups, 
\textit{Comm. in Algebra} {\bf 34} (2006), 
197--217. 


\bibitem{SV}
{\sc T.A. Springer} and {\sc F.D. Veldkamp}, 
\textit{Octonions, Jordan algebras and Exceptional Groups}, 
Springer-Verlag, 2000. 

\bibitem{Yi1959}
{\sc I.~Yokota}, 
Embeddings of projective spaces into 
elliptic projective Lie groups, 
{\it Proc. Japan Acad.} {\bf 35}-6 (1959), 
281--283. 

\bibitem{Yi1967} 
{\sc I.~Yokota}, 
Representation rings of group $G_2$, 
\textit{J. Fac. Sci. Shinshu Univ.} {\bf 2}-2 (1967), 
125--138.  


\bibitem{Yi1968} 
{\sc I.~Yokota}, 
Exceptional Lie group $F_4$ and its representation rings, 
\textit{J. Fac. Sci. Shinshu Univ.} {\bf 3}-1 (1968), 
35--60.  

\bibitem{Yi1975}
{\sc I.~Yokota}, 
On a non compact simple Lie group $F_{4, 1}$ of type $F_4$, 
{\it J. Fac. Sci. Shinshu Univ.} {\bf 10}-2 (1975), 71--80. 

\bibitem{Yi1977} 
{\sc I. Yokota}, 
Non-compact simple Lie group $F_{4, 2}$ of type $F_4$, 
\textit{J. Fac. Sci. Shinshu Univ.} {\bf 12} (1977), 53--64. 

\bibitem{Yi1990}
{\sc I.~Yokota}, 
Realizations of involutive automorphisms $\sigma$ and $G^\sigma$ 
of exceptional linear Lie groups $G$, Part \Roman{f1}, 
$G=G_2$. $F_4$ and $E_6$, 
\textit{Tsukuba J. Math.} {\bf 14}-1 (1990), 185--223.

\end{thebibliography}
\end{document}